\pdfoutput=1
\RequirePackage{ifpdf}
\ifpdf 
\documentclass[pdftex]{sigma}
\else
\documentclass{sigma}
\fi

\numberwithin{equation}{section}
\usepackage{mathrsfs}

\renewcommand\le{\leqslant}
\renewcommand\ge{\geqslant}

\newcommand{\Zset}{\mathbb Z}
\newcommand{\Rset}{\mathbb R}
\newcommand{\Cset}{\mathbb C}

\renewcommand{\Re}{\mathop{\rm Re}\nolimits}
\renewcommand{\Im}{\mathop{\rm Im}\nolimits}

\def\Ai{{\rm Ai}}
\def\Bi{{\rm Bi}}
\newcommand{\Beta}{{\rm B}}

\newcommand{\Lscr}{{\mathscr L}}

\newcommand{\lshad}{[\![}
\newcommand{\rshad}{]\!]}
\newcommand{\rme}{\mathrm{e}}
\newcommand{\rmi}{\mathrm{i}}
\newcommand{\rmd}{\mathrm{d}}

\newcommand{\dbltilde}[1]{\tilde{\raisebox{0pt}[0.85\height]{$\tilde{#1}$}}}

\newtheorem{Corollary}{Corollary}[section]

\begin{document}

\allowdisplaybreaks

\newcommand{\arXivNumber}{1707.05218}

\renewcommand{\PaperNumber}{042}

\FirstPageHeading

\ShortArticleName{Higher Derivatives of Airy Functions and of their Products}

\ArticleName{Higher Derivatives of Airy Functions\\ and of their Products}

\Author{Eugeny G.~ABRAMOCHKIN and Evgeniya V.~RAZUEVA}

\AuthorNameForHeading{E.G.~Abramochkin and E.V.~Razueva}

\Address{Coherent Optics Lab, Lebedev Physical Institute, Samara, 443011, Russia}
\Email{\href{mailto:ega@fian.smr.ru}{ega@fian.smr.ru}, \href{mailto:dev@fian.smr.ru}{dev@fian.smr.ru}}

\ArticleDates{Received October 13, 2017, in final form April 26, 2018; Published online May 05, 2018}

\Abstract{The problem of evaluation of higher derivatives of Airy functions in a closed form is investigated. General expressions for the polynomials which have arisen in explicit formulae for these derivatives are given in terms of particular values of Gegenbauer polynomials. Similar problem for products of Airy functions is solved in terms of terminating hypergeometric series.}


\Classification{33C10; 33C05; 33C20}

\section{Introduction}\label{sec:Intro}

Explicit formulae for higher derivatives is a usual part of investigation of special functions in mathematical physics \cite{BatemanI,BatemanII,BatemanIII, Szego,Watson, Whittaker}. A collection of these results can be found in~\cite{Brychkov, DLMF}.

For any solution $y(x)$ of a homogeneous differential equation of second order
\begin{gather*}
y''+p(x)y'+q(x)y=0
\end{gather*}
one can obtain that
\begin{gather*}
y''=-qy-py',\\
y'''=(pq-q')y+\big(p^2-p'-q\big)y',\\
y^{\rm IV}=\big({-}p^2q+pq'+2p'q+q^2-q''\big)y+\big({-}p^3+3pp'+2pq-p''-2q'\big)y',
\end{gather*}
etc. However, it would be quite difficult to get the explicit formula in general:
\begin{gather*}
y^{(n)}=P_n(x;p,q)y+Q_n(x; p,q)y',
\end{gather*}
because successful finding of coefficients $P_n(x;p,q)$ and $Q_n(x;p,q)$ in a closed form depends on many circumstances.

The main purpose of this paper is to obtain general formulae for $n$-th derivatives of Airy functions, $\Ai(x)$ and $\Bi(x)$. Both functions satisfy the same equation
\begin{gather}
y''=xy.\label{eq:AiryDE}
\end{gather}
Therefore,
\begin{gather}
\Ai^{(n)}(x)=P_n(x)\Ai(x)+Q_n(x)\Ai'(x),\nonumber\\
\Bi^{(n)}(x)=P_n(x)\Bi(x)+Q_n(x)\Bi'(x),\label{eq:DerivAi}
\end{gather}
where $P_n(x)$ and $Q_n(x)$ are some polynomials and the index $n$ corresponds to the derivative order but not the polynomials degree. By differentiating the first equation of (\ref{eq:DerivAi}),
\begin{gather*}
\Ai^{(n+1)}(x)=\bigl[P_n'(x)+xQ_n(x)\bigr]\Ai(x)+\bigl[P_n(x)+Q_n'(x)\bigr]\Ai'(x) \\
\hphantom{\Ai^{(n+1)}(x)}{}=P_{n+1}(x)\Ai(x)+Q_{n+1}(x)\Ai'(x),
\end{gather*}
we have two differential difference relations
\begin{gather}
P_{n+1}(x)=P_n'(x)+xQ_n(x),\qquad Q_{n+1}(x)=P_n(x)+Q_n'(x) \label{eq:DiffRecur_PQ}
\end{gather}
with initial conditions $P_0(x)=1$ and $Q_0(x)=0$ which help to determine $P_n(x)$ and $Q_n(x)$ for any natural $n$ (see the Table~\ref{table1}).

\begin{table}[t]\centering
\caption{First 16 polynomials $P_n(x)$ and $Q_n(x)$.}\label{table1}\vspace{1mm}

\begin{tabular}{rcc}
 $n$\qquad & $P_n(x)$ & $Q_n(x)$\\
 \hline
 $0$\qquad & $1$ & $0$\medskip\\
 $1$\qquad & $0$ & $1$\\
 $2$\qquad & $x$ & $0$\\
 $3$\qquad & $1$ & $x$\\
 $4$\qquad & $x^2$ & $2$\\
 $5$\qquad & $4x$ & $x^2$\medskip\\
 $6$\qquad & $x^3+4$ & $6x$\\
 $7$\qquad & $9x^2$ & $x^3+10$\\
 $8$\qquad & $x^4+28x$ & $12x^2$\\
 $9$\qquad & $16x^3+28$ & $x^4+52x$\\
$10$\qquad & $x^5+100x^2$ & $20x^3+80$\medskip\\
$11$\qquad & $25x^4+280x$ & $x^5+160x^2$\\
$12$\qquad & $x^6+260x^3+280$ & $30x^4+600x$\\
$13$\qquad & $36x^5+1380x^2$ & $x^6+380x^3+880$\\
$14$\qquad & $x^7+560x^4+3640x$ & $42x^5+2520x^2$\\
$15$\qquad & $49x^6+4760x^3+3640$ & $x^7+770x^4+8680x$
\end{tabular}
\end{table}

It is quite possible that the problem of evaluation of polynomials $P_n(x)$ and $Q_n(x)$ was formulated in XIX century when G.B.~Airy introduced the function which is now denoted $\Ai(x)$ (see~\cite{Airy1,Airy2}). However, as far as we know, the first solution of the problem has been published by Maurone and Phares in~1979 in terms of double finite sums containing factorials and gamma function ratio~\cite{Maurone}. We rewrite it using binomial coefficients and Pochhammer symbols
\begin{gather*}
P_{3m+\delta}(x)=\sum_{0\le k\le \lfloor\frac{m-k_1}2\rfloor}\frac{3^{m+m_1+k}x^{3k+\ell_1}}{(3k+\ell_1)!}\sum_{0\le\ell\le 3k+\ell_1} (-1)^\ell\binom{3k+\ell_1}{\ell}\left(\frac{1-\ell}{3}\right)_{m+m_1+k},\nonumber\\
Q_{3m+\delta}(x)=\sum_{0\le k\le \lfloor\frac{m-k_2}2\rfloor}\frac{3^{m+m_2+k}x^{3k+\ell_2}}{(3k+\ell_2)!}\sum_{0\le\ell\le 3k+\ell_2} (-1)^\ell\binom{3k+\ell_2}{\ell}\left(\frac{2-\ell}{3}\right)_{m+m_2+k},
\end{gather*}
where $\lfloor x\rfloor$ is the integral part of $x$, $\delta=0,1,2$, and
\begin{alignat*}{8}
& k_1=0, \quad && \ell_1=0, \quad && m_1=0,\quad && k_2=1, \quad && \ell_2=1, \quad && m_2=0,\quad && \mbox{if $\delta=0$},& \\
& k_1=1, \quad && \ell_1=2, \quad && m_1=1,\quad && k_2=0, \quad && \ell_2=0, \quad && m_2=0, \quad && \mbox{if $\delta=1$},& \\
& k_1=0, \quad && \ell_1=1, \quad && m_1=1,\quad && k_2=1, \quad && \ell_2=2, \quad && m_2=1,\quad && \mbox{if $\delta=2$}.&
\end{alignat*}

Recently, the problem of evaluation of $P_n(x)$ and $Q_n(x)$ has been investigated by Laurenzi in~\cite{Laurenzi} where the solution has been written in terms of Bell polynomials. More general problem of evaluation of higher derivatives of Bessel and Macdonald functions of arbitrary order has been solved by Brychkov in~\cite{Brychkov2013}. Corresponding polynomials have been written as finite sums containing products of terminated hypergeometric
series ${}_3F_2$ and ${}_2F_3$ but they look quite cumbersome for the case of~$\Ai(x)$, namely,
\begin{gather*}
P_n(x)=\frac{n!}{2x^{n-3}}\sum_{\lfloor\frac{n+2}{3}\rfloor\le k\le n}\frac{(-1)^{k+1}}{k!}\binom{3k}{n}\left(-\frac13\right)_{k}\\
\hphantom{P_n(x)=}{} \times{}_3F_2\left(\begin{matrix}
\frac{n}{3}-k,\frac{n+1}{3}-k,\frac{n+2}{3}-k \\
\frac13-k,\frac23-k \end{matrix}\,\Bigl|\,1\right)\cdot{}_2F_3\left(
\begin{matrix} 1-\frac{k}{2},\frac{3-k}{2} \\
2-k,\frac43-k,\frac53 \end{matrix}\,\Bigl|\,\frac{4x^3}{9}\right),\\
Q_n(x) =\frac{n!}{x^{n-1}}\sum_{\lfloor\frac{n+2}{3}\rfloor\le k\le n}
\frac{(-1)^k}{k!}\binom{3k}{n}\left(-\frac13\right)_{k}\\
\hphantom{Q_n(x) =} {}\times{}_3F_2\left(\begin{matrix}
\frac{n}{3}-k,\frac{n+1}{3}-k,\frac{n+2}{3}-k \\
\frac13-k,\frac23-k \end{matrix}\,\Bigl|\,1\right)\cdot{}_2F_3\left(
\begin{matrix} 1-\frac{k}{2},\frac{1-k}{2} \\
1-k,\frac43-k,\frac23 \end{matrix}\,\Bigl|\,\frac{4x^3}{9}\right),
\end{gather*}
where $n\ge 4$. However, much more cumbersome expression is produced by the on-line service {\it Wolfram Alpha} in response to a query ``$n$-th derivative of Airy function'' (see~\cite{Wolfram}).

This paper is organized as follows. In Section~\ref{sec:Ai_Deriv} we find simple expressions for $P_n(x)$ and~$Q_n(x)$ containing a special value of Gegenbauer polynomials. In Section~\ref{sec:HG_series} some corollaries of this result are considered. In Section~\ref{sec:Ai2_Deriv} we study higher derivatives of $\Ai^2(x)$, $\Ai(x)\Bi(x)$, and~$\Bi^2(x)$. In the last section we discuss some open problems connecting with the above results.

It is interesting that the evaluation problem for higher derivatives of $\Ai(x)$ and $\Ai^2(x)$ arise in physics for describing bound state solutions of the Schr\"odinger equation with a linear potential~\cite{Ferreira} and for quantum corrections of the Thomas--Fermi statistical model of atom~\cite{Schwinger}.

\section[Higher derivatives of $\Ai(x)$ and $\Bi(x)$]{Higher derivatives of $\boldsymbol{\Ai(x)}$ and $\boldsymbol{\Bi(x)}$}\label{sec:Ai_Deriv}

For evaluation of polynomials $P_n(x)$ and $Q_n(x)$ it is convenient to use {\it Airy atoms} $f(x)$ and $g(x)$ which are connected with Airy functions by the relations~\cite{DLMF}
\begin{gather*}
\Ai(x)=c_1f(x)-c_2g(x),\qquad \frac{\Bi(x)}{\sqrt 3}=c_1f(x)+c_2g(x),
\end{gather*}
where $c_1=3^{-2/3}/\Gamma\bigl(\tfrac23\bigr)$ and $c_2=3^{-1/3}/\Gamma\bigl(\tfrac13\bigr)$.

The reasons are the following. First of all, $f(x)$ and $g(x)$ are solutions of (\ref{eq:AiryDE}). As result,
\begin{gather}
f^{(n)}(x)=P_n(x)f(x)+Q_n(x)f'(x),\nonumber\\
g^{(n)}(x)=P_n(x)g(x)+Q_n(x)g'(x).\label{eq:AiryAtoms_PQ}
\end{gather}
Second, both Airy atoms have a simple hypergeometric representation:
\begin{gather*}
f(x)=\sum_{k=0}^\infty \left(\frac13\right)_k\frac{3^k x^{3k}}{(3k)!}={}_0F_1\left(\frac23\,\Bigl|\,\frac{x^3}{9}\right),\\
g(x)=\sum_{k=0}^\infty \left(\frac23\right)_k\frac{3^k x^{3k+1}}{(3k+1)!}=x\cdot{}_0F_1\left(\frac43\,\Bigl|\,\frac{x^3}{9}\right).
\end{gather*}
Then
\begin{gather*}
f^{(n)}(x)=\sum_{3k\ge n} \left(\frac13\right)_k \frac{3^k x^{3k-n}}{(3k-n)!},\qquad
g^{(n)}(x)=\sum_{3k+1\ge n} \left(\frac23\right)_k \frac{3^k x^{3k+1-n}}{(3k+1-n)!}.
\end{gather*}
And third, the Wronskian of $f(x)$ and $g(x)$ is equal to 1,
\begin{gather*}
W[f,g]=f(x)g'(x)-g(x)f'(x)\equiv 1,
\end{gather*}
that helps to simplify polynomials $P_n(x)$ and $Q_n(x)$ as solutions of the system (\ref{eq:AiryAtoms_PQ})
\begin{gather}
P_n(x)=g'(x)f^{(n)}(x)-f'(x)g^{(n)}(x),\nonumber\\
Q_n(x)=f(x)g^{(n)}(x)-g(x)f^{(n)}(x).\label{eq:PQ_fg}
\end{gather}
Substitution of the series expansions of Airy atoms and their derivatives into~(\ref{eq:PQ_fg}) yields
\begin{gather}
P_n(x)=\sum_{3m\ge n} \frac{3^m x^{3m-n}}{(3m-n)!}
\sum_k \left\{\binom{3m-n}{3k-n}-\binom{3m-n}{3k-1}\right\}
\left(\frac13\right)_k\left(\frac23\right)_{m-k},\label{eq:PQ_Series}\\
Q_n(x) =\sum_{3m\ge n-1} \frac{3^m x^{3m+1-n}}{(3m+1-n)!}
\sum_k \left\{ \binom{3m+1-n}{3k}-\binom{3m+1-n}{3k-n}\right\}
\left(\frac13\right)_k\left(\frac23\right)_{m-k},\nonumber
\end{gather}
where both series on $k$ are naturally terminating, i.e., the index of summation runs over all values that produce non-zero summands.

Let us define $\gamma_m$ by the formula
\begin{gather*}
\gamma_m(m_0,k_0)=\sum_k \binom{3m-m_0}{3k-k_0} \left(\frac13\right)_k\left(\frac23\right)_{m-k}.
\end{gather*}
Then (\ref{eq:PQ_Series}) can be rewritten as follows
\begin{gather}
P_n(x)=\sum_{3m\ge n}\frac{3^m x^{3m-n}}{(3m-n)!}\bigl[\gamma_m(n,n)-\gamma_m(n,1)\bigr],\nonumber\\
Q_n(x)=\sum_{3m\ge n-1}\frac{3^m x^{3m+1-n}}{(3m+1-n)!}\bigl[\gamma_m(n-1,0)-\gamma_m(n-1,n)\bigr].\label{eq:PQ_Series2}
\end{gather}

Reducing the product of Pochhammer symbols to Euler's beta function
\begin{gather*}
\left(\frac13\right)_k\left(\frac23\right)_{m-k}=\frac{\Gamma\bigl(k+\tfrac13\bigr)\Gamma\bigl(m-k+\tfrac23\bigr)}
{\Gamma\bigl(\tfrac13\bigr)\Gamma\bigl(\tfrac23\bigr)}=\frac{m!\sqrt 3}{2\pi} \Beta\left(k+\frac13,m-k+\frac23\right),
\end{gather*}
and applying one of the function integral definitions
\begin{gather*}
\Beta(x,y)=\int_0^\infty \!\frac{\tau^{x-1}}{(1+\tau)^{x+y}}\,\rmd\tau, \qquad \Re x>0, \qquad \Re y>0,
\end{gather*}
we have
\begin{gather}
\gamma_m(m_0,k_0)=\frac{m!\sqrt 3}{2\pi}\int_0^\infty \sum_k\binom{3m-m_0}{3k-k_0}\tau^k\cdot\frac{\tau^{-2/3}}{(1+\tau)^{m+1}}\,\rmd\tau\nonumber\\
\hphantom{\gamma_m(m_0,k_0)} {}=\frac{m!\sqrt 3}{2\pi}\int_0^\infty \sum_k\binom{3m-m_0}{3k-k_0} t^{3k}\cdot\frac{3}{(1+t^3)^{m+1}}\,\rmd t\nonumber\\
\hphantom{\gamma_m(m_0,k_0)} {}=\frac{m!\sqrt 3}{2\pi}\int_0^\infty \sum_k \binom{3m-m_0}{k}
\bigl[1+\omega^{k+k_0}+\bar\omega^{k+k_0}\bigr]\frac{t^{k+k_0}}{(1+t^3)^{m+1}}\,\rmd t\nonumber\\
\hphantom{\gamma_m(m_0,k_0)} {}=\frac{m!\sqrt 3}{2\pi}\sum_{0\le j\le 2} \omega^{jk_0}\int_0^\infty
\frac{t^{k_0}(1+\omega^j t)^{3m-m_0}}{(1+t^3)^{m+1}}\,\rmd t.\label{eq:H1}
\end{gather}
Here, $\omega=\rme^{2\pi\rmi/3}$ and $\bar\omega=\rme^{-2\pi\rmi/3}$ are the cubic roots of unity, and an overline means complex conjugation.

Changing the variable in the last integral, $t\to 1/t$, we get another expression for $\gamma_m$:
\begin{gather}
\gamma_m(m_0,k_0)=\frac{m!\sqrt 3}{2\pi}\sum_{0\le j\le 2} \omega^{j(k_0-m_0)}\int_0^\infty \frac{t^{m_0-k_0+1}(1+\bar\omega^j t)^{3m-m_0}}{(1+t^3)^{m+1}}\,\rmd t.\label{eq:H2}
\end{gather}
Then the differences in square brackets in (\ref{eq:PQ_Series2}) can be written as
\begin{gather*}
\gamma_m(n,n)\Bigl|_{(\ref{eq:H2})}{}-\gamma_m(n,1)\Bigl|_{(\ref{eq:H1})}
=\frac{m!\sqrt 3}{2\pi}\cdot 2\Re\biggl\{(1-\omega)\int_0^\infty t
\frac{(1+\omega t)^{3m-n}}{(1+t^3)^{m+1}}\,\rmd t\biggr\},\\
\gamma_m(n-1,0)\Bigl|_{(\ref{eq:H1})}{}-\gamma_m(n-1,n)\Bigl|_{(\ref{eq:H2})}
=\frac{m!\sqrt 3}{2\pi}\cdot 2\Re\biggl\{(1-\bar\omega)\int_0^\infty
\frac{(1+\omega t)^{3m+1-n}}{(1+t^3)^{m+1}}\,\rmd t\biggr\}.
\end{gather*}
Expansions (\ref{eq:PQ_Series2}) show that exponents of the factor $(1+\omega t)$ in both integrals are nonnegative.

The integral for $\gamma_m(n-1,0)-\gamma_m(n-1,n)$ looks a little simpler than for $\gamma_m(n,n)-\gamma_m(n,1)$, so we begin with it. Integrating the analytic function $(1+z)^{3m+1-n}/(1+z^3)^{m+1}$ over the contour $[0,R]\cup \{ R\rme^{\rmi\phi},\,\phi\in[0,2\pi/3]\}\cup[R\omega,0]$ and taking $R\to\infty$, we get
\begin{gather*}
\int_0^\infty \frac{(1+t)^{3m+1-n}}{(1+t^3)^{m+1}}\,\rmd t +\int_\infty^0 \frac{(1+\omega t)^{3m+1-n}}{(1+t^3)^{m+1}}\omega\,\rmd t
=2\pi\rmi\mathop{\rm res}_{z=\rme^{\pi\rmi/3}} \frac{(1+z)^{3m+1-n}}{(1+z^3)^{m+1}}.
\end{gather*}
Then the integral of our interest can be written as
\begin{gather*}
\int_0^\infty \frac{(1+\omega t)^{3m+1-n}}{(1+t^3)^{m+1}}\,\rmd t
=\bar\omega \int_0^\infty \frac{(1+t)^{3m+1-n}}{(1+t^3)^{m+1}}\,\rmd t
-\frac{2\pi\rmi\bar\omega}{m!}\frac{\rmd^m}{\rmd t^m}
\left(\frac{(1+t)^{2m-n}}{(t-\rme^{-\pi\rmi/3})^{m+1}}\right) \biggl|_{t=\rme^{\pi\rmi/3}}.
\end{gather*}
It is worth noting that the first term on the right side makes a zero contribution to the difference $\gamma_m(n-1,0)-\gamma_m(n-1,n)$ because $\Re\{(1-\bar\omega)\bar\omega\}=\Re(-\rmi\sqrt 3)=0$. Thus,
\begin{gather}
\gamma_m(n-1,0)-\gamma_m(n-1,n) =-6\Re\frac{\rmd^m}{\rmd t^m} \left(\frac{(1+t)^{2m-n}}{(t-\rme^{-\pi\rmi/3})^{m+1}}\right)\biggl|_{t=\rme^{\pi\rmi/3}}\nonumber\\
\hphantom{\gamma_m(n-1,0)-\gamma_m(n-1,n)} {}=-\frac{2}{(\sqrt{3})^{n-1}}\Re\frac{\rmd^m}{\rmd t^m}\left(\frac{(t+\rme^{\pi\rmi/6})^{2m-n}}{(t+\rmi)^{m+1}}\right)\biggl|_{t=0}.\label{eq:H1a}
\end{gather}
Using the generating function, we prove that the right part of (\ref{eq:H1a}) vanishes for $n=0,1,\ldots,2m$. Let $s\in\Rset$. Then
\begin{gather*}
F(s) =\sum_{n=0}^{2m} \binom{2m}{n}\frac{\rmd^m}{\rmd t^m}\left(\frac{(t+\rme^{\pi\rmi/6})^{2m-n}}{(t+\rmi)^{m+1}}\right)\biggl|_{t=0}{}\cdot s^n\\
\hphantom{F(s)} {}=\frac{\rmd^m}{\rmd t^m} \left(\frac{(t+\rme^{\pi\rmi/6}+s)^{2m}}{(t+\rmi)^{m+1}}\right)\biggl|_{t=0}=\frac{m!}{2\pi\rmi}\oint_{|z|=\epsilon}\frac{(z+\rme^{\pi\rmi/6}+s)^{2m}}{(z+\rmi)^{m+1}z^{m+1}}\,\rmd z.
\end{gather*}
Changing the variable $\zeta=(z+a)/(1+bz)$, where parameters $a$ and $b$ will be chosen later, one gets
\begin{gather*}
F(s)=\frac{m!}{2\pi\rmi}(1-ab)\oint_{|\zeta-a|=\epsilon} \frac{(\zeta[1-b\rme^{\pi\rmi/6}-bs]+\rme^{\pi\rmi/6}+s-a)^{2m}}{(\zeta[1-\rmi b]+\rmi-a)^{m+1}(\zeta-a)^{m+1}}\,\rmd\zeta.
\end{gather*}
With the aim to simplify the integrand, we set $b=-\rmi$ and $a=\rme^{\pi\rmi/6}+s$. Then
\begin{gather*}
1-b\rme^{\pi\rmi/6}-bs=1-ab=\rmi\bar a,\qquad \rmi-a=-\bar a,
\end{gather*}
and
\begin{gather*}
F(s)=\frac{m!}{2\pi\rmi} \frac{(\rmi\bar a)^{2m+1}}{(-\bar a)^{m+1}} \oint_{|\zeta-a|=\epsilon} \frac{\zeta^{2m}}{(\zeta-a)^{m+1}}\,\rmd\zeta=-\rmi\bar a^m\frac{\rmd^m}{\rmd t^m} t^{2m}\biggl|_{t=a}=-\rmi \frac{(2m)!}{m!} |a|^{2m}.
\end{gather*}
Therefore, $\Re F(s)=0$, and the expansion of $Q_n(x)$ in~(\ref{eq:PQ_Series2}) can be simplified by throwing out the terms with $n<2m+1$. As result, the function
\begin{gather}
Q_n(x)=\sum_{\frac{n-1}{3}\le m\le\frac{n-1}{2}} \frac{x^{3m+1-n}}{(3m+1-n)!}\nonumber\\
\hphantom{Q_n(x)=\sum_{\frac{n-1}{3}\le m\le\frac{n-1}{2}} }{} \times\frac{-2}{(\sqrt{3})^{n-2m-1}}
\Re\frac{\rmd^m}{\rmd t^m}\left(\frac{1}{(t+\rmi)^{m+1}(t+\rme^{\pi\rmi/6})^{n-2m}}\right)\!\biggl|_{t=0}\label{eq:Q_fg3}
\end{gather}
is actually a polynomial.

The double inequality $\frac{n-1}{3}\le m\le\frac{n-1}{2}$ is quite restrictive for the summation index. In particular, it leads immediately to zero polynomials when $n$ is equal to $0$ or~$2$.

Now we need to evaluate the derivative in (\ref{eq:Q_fg3}). As before, we use the generating function approach. Let
\begin{gather*}
g_{m,n}=\frac{1}{m!} \frac{\rmd^m}{\rmd t^m} \left(\frac{1}{(t+\rmi)^{m+1}(t+\rme^{\pi\rmi/6})^{n+1}}\right)\biggl|_{t=0}, \qquad m,n\ge 0
\end{gather*}
and $s\in\Rset$. Then
\begin{gather*}
G(s)=\sum_{m=0}^\infty g_{m,n}s^m =\frac{1}{2\pi\rmi}\oint_{|z|=\epsilon}\sum_{m=0}^\infty \,\frac{s^m}{z^{m+1}(z+\rmi)^{m+1}}\cdot\frac{\rmd z}{(z+\rme^{\pi\rmi/6})^{n+1}}\\
\hphantom{G(s)=\sum_{m=0}^\infty g_{m,n}s^m}{}=\frac{1}{2\pi\rmi}\oint_{|z|=\epsilon} \frac{\rmd z}{(z^2+\rmi z-s)(z+\rme^{\pi\rmi/6})^{n+1}},
\end{gather*}
assuming that $\epsilon$ is a small positive number and $|s|<\min_{|z|=\epsilon} |z(z+\rmi)|=\epsilon(1-\epsilon)$ for the geometric series convergence. Roots of the equation $z^2+\rmi z-s=0$ are $z_{\pm}=-\rmi(1\pm\sqrt{1-4s}\,)/2$, and $z_{-}$ is the only pole of the integrand within the contour $|z|=\epsilon$. Therefore,
\begin{gather*}
G(s)=\mathop{\rm res}_{z=z_{-}}\frac{1}{(z^2+\rmi z-s)(z+\rme^{\pi\rmi/6})^{n+1}}=\frac{2^{n+1}}{\rmi\sqrt{1-4s}(\sqrt{3}+\rmi\sqrt{1-4s})^{n+1}}, \\
g_{m,n}=\lshad s^m\rshad G(s)=\left(\frac{2}{\sqrt 3}\right)^{n+1}\lshad s^m\rshad\left\{\frac{1}{\rmi\sqrt{1-4s}}\sum_{k=0}^\infty
\binom{n+k}{n} \left(-\frac{\rmi\sqrt{1-4s}}{\sqrt 3}\right)^{k}\right\}.
\end{gather*}
Here, we use the notation proposed in \cite{GKP}. Namely, if $A(z)$ is any power series $\sum_k a_kz^k$, then $\lshad z^k\rshad A(z)$ denotes the coefficient of $z^k$ in $A(z)$. In our view, this notation is more convenient to manipulate power series than usual analytic description, $\lshad z^k\rshad A(z)=A^{(k)}(0)/k!$.

Returning to (\ref{eq:Q_fg3}) and renaming the coefficients for brevity, we have
\begin{gather}
\tilde g_{m,n} =\frac{-2}{(\sqrt{3})^n}\Re g_{m,n} =\frac{2^{n+2}}{3^{n+1}}\,\lshad s^m\rshad \sum_{k=0}^\infty \binom{n+2k+1}{n}\left(-\frac{1-4s}{3}\right)^{k}
\nonumber\\
\hphantom{\tilde g_{m,n}} {}=\frac{2^{n+2}}{3^{n+1}}\sum_{k=m}^\infty \binom{n+2k+1}{n}\left(-\frac13\right)^{k}\cdot\binom{k}{m}(-4)^m\label{eq:BinomialSeries}\\
\hphantom{\tilde g_{m,n}} {}=\frac{2^{n+2m+2}}{3^{n+m+1}}\binom{n+2m+1}{n}\cdot {}_2F_1 \left(\frac{n+2m+3}{2},\frac{n+2m+2}{2}; \frac{2m+3}{2}\,\Bigl|\,-\frac13\right).\nonumber
\end{gather}

Using Euler's transformation ${}_2F_1(a,b;\,c\,|\,z)=(1-z)^{c-a-b}\cdot{}_2F_1(c-a,c-b;\,c\,|\,z)$ and one of hypergeometric representations of Gegenbauer's polynomials \cite[equation (7.3.1.202)]{Prudnikov3},
\begin{gather*}
{}_2F_1 \left(-\frac{n}{2},-\frac{n-1}{2}; \lambda+\frac12\,\Bigl|\, 1-z^2\right)=\frac{n!\,z^n}{(2\lambda)_n}\,C_n^\lambda\left(\frac{1}{z}\right),
\end{gather*}
we can simplify the coefficients to
\begin{gather}
\tilde g_{m,n} =\frac{1}{2^n}\binom{n+2m+1}{n}\cdot {}_2F_1\left(-\frac{n}{2},-\frac{n-1}{2};m+\frac32\,\Bigl|\,-\frac13\right)\label{eq:HyperG_Q}\\
\hphantom{\tilde g_{m,n}}{}=\frac{1}{(\sqrt{3})^n} C_n^{m+1}\left(\frac{\sqrt 3}{2}\right) =\lshad t^n\rshad\frac{1}{\bigl(1-t+\tfrac13 t^2\bigr)^{m+1}}\label{eq:Q_Coeffs}
\end{gather}
and obtain required formulae for both polynomial families, first for $Q_n(x)$, then for $P_n(x)=Q_{n+1}(x)-Q_n'(x)$. We write them in the same style, shifting the index in $Q_n(x)$
\begin{gather}
Q_{n+1}(x)=\sum_{\frac{n}{3}\le m\le\frac{n}{2}} \tilde g_{m,n-2m}\frac{m!x^{3m-n}}{(3m-n)!},\label{eq:Solution_Q}\\
P_n(x)=\sum_{\frac{n}{3}\le m\le\frac{n}{2}}\{\tilde g_{m,n-2m}-\tilde g_{m,n-2m-1}\}\frac{m! x^{3m-n}}{(3m-n)!}.\label{eq:Solution_P}
\end{gather}
The last expression in (\ref{eq:Q_Coeffs}) follows from the familiar generating function for Gegenbauer polynomials~\cite{DLMF}
\begin{gather*}
\sum_{n=0}^\infty C_n^\lambda(x)t^n=\frac{1}{(1-2xt+t^2)^\lambda}, \qquad |t|<1,\qquad |x|<1.
\end{gather*}
Of course, all the coefficients $\tilde g_{m,n-2m}$ and $\tilde g_{m,n-2m} -\tilde g_{m,n-2m-1}$ are positive but this is easier to deduce from (\ref{eq:DiffRecur_PQ}) than from (\ref{eq:Q_Coeffs}). It is worth noting, that the term $\tilde g_{m,n-2m-1}$ in (\ref{eq:Solution_P}) vanishes for the case $n=2m$ since (\ref{eq:HyperG_Q}) holds.

Left and right inequalities on the summation index in (\ref{eq:Solution_Q}) and (\ref{eq:Solution_P}) reveal the different structures of analytic expressions of the polynomials for small and large values of~$x$. Namely, for $x\sim 0$ the expressions depend on $n\,({\rm mod}\,3)$:
\begin{gather}
P_{3n}(x)=3^n\left(\frac13\right)_{n}+3^n\left\{(n+1)\left(\frac13\right)_{n}-\left(\frac23\right)_{n}\right\}x^3+\cdots,\nonumber\\
P_{3n+1}(x)=\frac{3^n}{2}\left\{\left(\frac43\right)_{n}-\left(\frac23\right)_{n}\right\}x^2+\frac{3^n}{40}\left\{(n+8)\left(\frac43\right)_{n}-(10n+8)\left(\frac23\right)_{n}\right\}x^5+\cdots,\nonumber\\
P_{3n+2}(x)=3^n\left(\frac43\right)_{n}x+\frac{3^n}{8}\left\{(n+4)\left(\frac43\right)_{n}-4\left(\frac53\right)_{n}\right\}x^4+\cdots,\nonumber\\
Q_{3n}(x)=3^n\left\{\left(\frac23\right)_{n}-\left(\frac13\right)_{n}\right\}x+\frac{3^n}{8}\left\{(n+2)\left(\frac23\right)_{n}-(4n+2)\left(\frac13\right)_{n}\right\}x^4+\cdots,\nonumber\\
Q_{3n+1}(x)=3^n\left(\frac23\right)_{n}+\frac{3^n}{2}\left\{(n+1)\left(\frac23\right)_{n}-\left(\frac43\right)_{n}\right\}x^3+\cdots,\label{eq:PQ1}\\
Q_{3n+2}(x)=3^n\left\{\left(\frac53\right)_{n}\!-\left(\frac43\right)_{n}\right\}x^2+\frac{3^n}{40}\left\{(2n+10)\left(\frac53\right)_{n}\!-(5n+10)\left(\frac43\right)_{n}\right\}x^5+\cdots,\nonumber
\end{gather}
while for $x\to\infty$ they depend on $n$ $({\rm mod}~2)$:
\begin{gather}
P_{2n}(x)=x^n+\binom{n}{3}(3n-5)x^{n-3} +10\binom{n}{6}(3n^2-15n+10)x^{n-6}+\cdots,\nonumber\\
P_{2n+1}(x) =n^2x^{n-1}+4\binom{n}{4}(n^2 - 2n - 1)x^{n-4} +14\binom{n}{7}(3n^3 - 17n^2+8n+8)x^{n-7}+\cdots,\nonumber\\
Q_{2n}(x) =2\binom{n}{2}x^{n-2}+20\binom{n}{5}(n-1)x^{n-5} +112\binom{n}{8}(3n-2)(n-3)x^{n-8}+\cdots,\nonumber\\
Q_{2n+1}(x) =x^n+\binom{n}{3}(3n+1)x^{n-3} +10\binom{n}{6}(3n^2-3n-2)x^{n-6}+\cdots. \label{eq:PQ2}
\end{gather}
In addition, it may be worth noting that the methods of finding these expansions are essentially different: (\ref{eq:PQ2}) follows directly from (\ref{eq:Solution_Q}) and (\ref{eq:Solution_P}), while for proving (\ref{eq:PQ1}) the easiest way is to apply (\ref{eq:PQ_Series}).

\section{Hypergeometric series and difference equations}\label{sec:HG_series}

The initial statement of the problem (\ref{eq:DerivAi}) and our final expressions for polynomials $P_n(x)$ and~$Q_n(x)$ lead to some corollaries whose relation to the Airy functions does not look evident a~priori. Below we consider two of them. The first helps to find the Gauss hypergeometric function special values. The second is connected with difference equations of third and fourth orders.

\subsection{Some special values of the hypergeometric function of Gauss}

Power series expansions of polynomials $P_n(x)$ and $Q_n(x)$ help to find values of the Gauss hypergeometric function in some special cases. The simplest of them is obtained equating the terms $Q_{3n+1}(0)$ in~(\ref{eq:Solution_Q}) and~(\ref{eq:PQ1}):
\begin{gather*}
Q_{3n+1}(0)=3^n\left(\frac23\right)_{n}=n!\tilde g_{n,n}=\frac{n!}{2^n}\binom{3n+1}{n}\cdot{}_2F_1\left(-\frac{n}{2},-\frac{n-1}{2};n+\frac32\,\Bigl|\,-\frac13\right).
\end{gather*}
Then
\begin{gather*}
{}_2F_1\left(-\frac{n}{2},-\frac{n-1}{2};\,n+\frac32\,\Bigl|\,-\frac13\right)
=6^n\left(\frac23\right)_{n}\cdot\frac{(2n+1)!}{(3n+1)!}=\left(\frac89\right)^{n}\cdot\frac{\bigl(\tfrac32\bigr)_n}{\bigl(\tfrac43\bigr)_n}.
\end{gather*}
This formula is not new (see Exercise~15 in \cite[Chapter~14]{Whittaker}) and can be extended to real or complex $a$'s by replacing $n$ by $-2a$ \cite[equation~(7.3.9.15)]{Prudnikov3}:
\begin{gather}
{}_2F_1\left(a,a+\frac12;\frac32-2a\,\Bigl|\,-\frac13\right)=\frac{\Gamma\bigl(\tfrac23-2a\bigr)\Gamma(2-4a)}{6^{2a}\Gamma\bigl(\tfrac23\bigr)\Gamma(2-6a)} =\left(\frac98\right)^{2a} \cdot \frac{2\Gamma\bigl(\tfrac32-2a\bigr)\Gamma\bigl(\tfrac43\bigr)} {\sqrt{\pi} \Gamma\bigl(\tfrac43-2a\bigr)}.\label{eq:HyperG_2F1a}
\end{gather}
As is typical in the theory of hypergeometric functions, the extension from~$n$ to $a$ is valid due to {\it Carlson's theorem} \cite[Section~5.8.1]{Titchmarsh_TheorFunc_ENG}: {\it If $f(z)$ is regular and of
the form ${\cal O}(\rme^{c|z|})$, where $c<\pi$, for $\Re z\ge 0$, and $f(z)=0$ for $z=0,1,2,\ldots$, then $f(z)=0$ identically.} (See \cite{Bailey, Borwein, Ebisu} for many examples of application of Carlson's theorem to
hypergeometric series.)

Next two hypergeometric function special values are related to coefficients of $Q_{3n+3}(x)$ and $Q_{3n+2}(x)$:
\begin{gather*}
\lshad x\rshad Q_{3n+3}(x) =3^{n+1}\left\{\left(\frac23\right)_{n+1} -\left(\frac13\right)_{n+1}\right\}=(n+1)! \tilde g_{n+1,n}\\
\hphantom{\lshad x\rshad Q_{3n+3}(x)}{}=\frac{(n+1)!}{2^n}\binom{3n+3}{n}\cdot{}_2F_1\left(-\frac{n}{2}, -\frac{n-1}{2};n+\frac52\,\Bigl|\,-\frac13\right),\\
\lshad x^2\rshad Q_{3n+2}(x)=3^n\left\{\left(\frac53\right)_{n}-\left(\frac43\right)_{n}\right\}=\frac{(n+1)!}{2}\tilde g_{n+1,n-1}\\
\hphantom{\lshad x^2\rshad Q_{3n+2}(x)}{}=\frac{(n+1)!}{2^n}\binom{3n+2}{n-1}\cdot{}_2F_1\left(-\frac{n-1}{2}, -\frac{n-2}{2};n+\frac52\,\Bigl|\,-\frac13\right).
\end{gather*}
By replacing $n$ by $-2a$ in the first relation and by $1-2a$ in the second, we get
\begin{gather*}
{}_2F_1\left(a,a+\frac12; \frac52-2a\,\Bigl|\,-\frac13\right) =\frac{6^{-2a}}{1-2a}\left\{2\frac{\Gamma\bigl(\tfrac53-2a\bigr)}{\Gamma\bigl(\tfrac53\bigr)}
-\frac{\Gamma\bigl(\tfrac43-2a\bigr)}{\Gamma\bigl(\tfrac43\bigr)}\right\} \frac{\Gamma(4-4a)}{\Gamma(4-6a)},\\ 
{}_2F_1\left(a,a+\frac12;\frac72-2a\,\Bigl|\,-\frac13\right)=\frac{6^{1-2a}}{(1-2a)(2-2a)}
\left\{\frac{\Gamma\bigl(\tfrac83-2a\bigr)}{\Gamma\bigl(\tfrac53\bigr)}
-\frac{\Gamma\bigl(\tfrac73-2a\bigr)}{\Gamma\bigl(\tfrac43\bigr)}\right\} \frac{\Gamma(6-4a)}{\Gamma(6-6a)}.
\end{gather*}

In order to use the coefficients of $P_n(x)$, we need to do some preliminary work. Substitution of (\ref{eq:BinomialSeries}) in both terms of the difference $\tilde g_{m,n}-\tilde g_{m,n-1}$ and the binomial identity
\begin{gather*}
\binom{N}{n-1}+\binom{N}{n}=\binom{N+1}{n}
\end{gather*}
lead to the relation
\begin{gather*}
\tilde g_{m,n}-\tilde g_{m,n-1} =\frac{3}{2^{n+1}}\binom{n+2m}{n}\cdot{}_2F_1\left(-\frac{n}{2},-\frac{n+1}{2}; m+\frac12\,\Bigl|\,-\frac13\right)\\
\hphantom{\tilde g_{m,n}-\tilde g_{m,n-1} =} {}-\frac{1}{2^{n+1}}\binom{n+2m+1}{n}\cdot{}_2F_1\left(-\frac{n-1}{2},-\frac{n}{2}; m+\frac32\,\Bigl|\,-\frac13\right).
\end{gather*}
Then, by considering the first nonzero coefficients of $P_{3n}(x)$ and
$P_{3n+2}(x)$,
\begin{gather*}
P_{3n}(0)=3^n\left(\frac13\right)_{n}=n!\{\tilde g_{n,n}-\tilde g_{n,n-1}\},\\
\lshad x\rshad P_{3n+2}(x) =3^n\left(\frac43\right)_{n}=(n+1)!\{\tilde g_{n+1,n}-\tilde g_{n+1,n-1}\},
\end{gather*}
we get the following special values of the Gauss hypergeometric function:
\begin{gather*}
{}_2F_1\left(a,a+\frac12; -\frac12-2a\,\Bigl|\,-\frac13\right) =\frac{6^{-2a}}{3}
\left\{\frac{\Gamma\bigl(\tfrac13-2a\bigr)}{\Gamma\bigl(\tfrac13\bigr)} +\frac{1+3a}{1+6a}\cdot
\frac{\Gamma\bigl(\tfrac23-2a\bigr)}{\Gamma\bigl(\tfrac23\bigr)}\right\} \frac{\Gamma(-1-4a)}{\Gamma(-1-6a)},\\
{}_2F_1\left(a,a+\frac12; \frac12-2a\,\Bigl|\,-\frac13\right) =\frac{6^{-2a}}{2}
\left\{\frac{\Gamma\bigl(\tfrac13-2a\bigr)}{\Gamma\bigl(\tfrac13\bigr)}
+\frac{\Gamma\bigl(\tfrac23-2a\bigr)}{\Gamma\bigl(\tfrac23\bigr)}\right\} \frac{\Gamma(1-4a)}{\Gamma(1-6a)}.
\end{gather*}

It is quite clear that this approach provides a way for evaluation of
\begin{gather*}
{}_2F_1\left(a,a+\frac12; n+\frac12-2a\,\Bigl|\,-\frac13\right)
\end{gather*}
in a closed form for any $n\in\Zset$, while, of course, the expressions will become more and more cumbersome as $|n|$ increases.

\subsection{Difference equations of third and fourth orders}

Relations (\ref{eq:PQ_fg}) help to find generating functions for polynomials $P_n(x)$ and $Q_n(x)$:
\begin{gather*}
{\bf P}(x,t)=\sum_{n=0}^\infty P_n(x)\frac{t^n}{n!}=g'(x)f(x+t)-f'(x)g(x+t),\\
{\bf Q}(x,t)=\sum_{n=0}^\infty Q_n(x)\frac{t^n}{n!}=f(x)g(x+t)-g(x)f(x+t).
\end{gather*}
Since $f(x)$ and $g(x)$ satisfy the equation (\ref{eq:AiryDE}), then ${\bf P}(x,t)$ and ${\bf Q}(x,t)$ being the functions on~$t$ are solutions of the differential equation
\begin{gather*}
\left\{\frac{\rmd^2}{\rmd t^2}-(x+t)\right\}{\bf Y}(x,t)=0.
\end{gather*}
If a solution of this equation is expanded in a power series
\begin{gather*}
{\bf Y}(x,t)=\sum_{n=0}^\infty Y_n(x)\frac{t^n}{n!},
\end{gather*}
then the series substitution into the equation leads to a recurrence relation for its coefficients, functions $Y_n(x)$:
\begin{gather}
Y_{n+3}(x)=xY_{n+1}(x)+(n+1)Y_n(x).\label{eq:Recur_PQZ}
\end{gather}
As result, polynomials $P_n(x)$ and $Q_n(x)$ are two of three solutions of~(\ref{eq:Recur_PQZ}), corresponding to different initial conditions:
\begin{alignat*}{5}
& P_n(x)\colon \quad && P_0(x)=1, \qquad && P_1(x)=0, \qquad && P_2(x)=x,& \\
& Q_n(x)\colon \quad && Q_0(x)=0, \qquad && Q_1(x)=1, \qquad && Q_2(x)=0.&
\end{alignat*}

It is not evident how to find the third independent solution of~(\ref{eq:Recur_PQZ}), say
\begin{gather*}
Z_n(x)\colon \quad Z_0(x)=0, \qquad Z_1(x)=0, \qquad Z_2(x)=1,
\end{gather*}
and how this relates with Airy functions' derivatives.

Investigation of the polynomials $Z_n(x)$ shows that
\begin{gather}
Z_{n+2}(x)=\sum_{\frac{n}{3}\le m\le\frac{n}{2}} \lambda_{m,n}\frac{m!x^{3m-n}}{(3m-n)!},\label{eq:Pol_Z}
\end{gather}
where coefficients $\lambda_{m,n}$ satisfy the relation
\begin{gather*}
m\lambda_{m,n}=(3m-n)\lambda_{m-1,n-2}+n\lambda_{m-1,n-3},\qquad m\ge 1.
\end{gather*}
In particular, for large values of $x$ we have
\begin{gather*}
Z_{2n+2}(x)=x^n+(n-1)(n-2)\frac{3n^2+7n+6}{6} x^{n-3}+\cdots,\\
Z_{2n+3}(x)=n(n+2)x^{n-1}+\binom{n-1}{3}(n+2)(n^2+2n+3)x^{n-4}+\cdots ,
\end{gather*}
while for $x\sim 0$ the first nonzero coefficients are
\begin{gather*}
\lshad x^2\rshad Z_{3n}(x)=\frac{3^n}{2}\left\{n!-2\left(\frac23\right)_{n}+\left(\frac13\right)_{n}\right\},\\
\lshad x\rshad Z_{3n+1}(x)=3^n\left\{n!-\left(\frac23\right)_{n}\right\}, \qquad Z_{3n+2}(0)=3^n n!.
\end{gather*}
Nevertheless, we could not reveal an analytic source of polynomials $Z_n(x)$ in order to find their expansions by applying the calculus methods as above.

We mention also two difference equations for the coefficients of polynomials $P_n(x)$ and $Q_n(x)$. These equations can be found considering the Laplace transform of the shifted Airy function,
\begin{gather*}
\Lscr(p,x)=\int_0^\infty \exp(-pt)\Ai(t+x)\,\rmd t.
\end{gather*}
We will use a formal power series approach for simplicity, while there is, of course, a more rigorous but tedious way based on an asymptotic behaviour of $\Ai(t)$ for large $t$'s. Assuming that $p\to +\infty$, the asymptotic expansion of $\Lscr(p,x)$ can be written in two ways. The first follows from~(\ref{eq:DerivAi})
\begin{gather}
\Lscr(p,x) =\sum_{n=0}^\infty \frac{\Ai^{(n)}(x)}{n!}\int_0^\infty \exp(-pt)t^n\,\rmd t=\sum_{n=0}^\infty \frac{P_n(x)\Ai(x)+Q_n(x)\Ai'(x)}{p^{n+1}}.\label{eq:AsympSeries1}
\end{gather}
On the other hand, the function $\Lscr(p,x)$ satisfies a differential equation of first order
\begin{gather}
\left\{\frac{\rmd}{\rmd p}+\big(p^2-x\big)\right\}\Lscr(p,x)=\Ai(x)p+\Ai'(x). \label{eq:DiffEq_LaplaceAi}
\end{gather}
With initial condition
\begin{gather*}
\Lscr(0,x)=\int_x^\infty \Ai(t)\,\rmd t=\Ai_1(x),
\end{gather*}
where the last notation is due to Aspnes \cite{Aspnes1}, the solution is
\begin{gather*}
\Lscr(p,x)=\exp\left(xp-\frac{p^3}{3}\right)\left\{\Ai_1(x) +\int_0^p \bigl(\Ai(x)t+\Ai'(x)\bigr)\exp\left(-xt+\frac{t^3}{3}\right)\rmd t\right\}.
\end{gather*}
Integration by parts yields
\begin{gather*}
\Lscr(p,x)=\Ai(x)\left\{\frac{p}{p^2-x}+\frac{1}{(p^2-x)^2}+\frac{2x}{(p^2-x)^3}+\frac{4p}{(p^2-x)^4}+\cdots\right\}\\
\hphantom{\Lscr(p,x)=}{}+\Ai'(x)\left\{\frac{1}{p^2-x}+\frac{2p}{(p^2-x)^3} +\frac{10}{(p^2-x)^4}+\ldots\right\}+{\cal O}\left(\frac{1}{p^\infty}\right).
\end{gather*}
Then, in general, $\Lscr(p,x)$ has the form
\begin{gather}
\Lscr(p,x)=\Ai(x)\sum_{k=0}^\infty \frac{\mu_k p+\nu_k}{(p^2-x)^{k+1}}+\Ai'(x)\sum_{k=0}^\infty \frac{\tilde\mu_k p+\tilde\nu_k}{(p^2-x)^{k+1}}\label{eq:AsympSeries2}
\end{gather}
with coefficients $\mu_k$, $\nu_k$, $\tilde\mu_k$, $\tilde\nu_k$ depending on~$x$.

Substituting (\ref{eq:AsympSeries2}) in (\ref{eq:DiffEq_LaplaceAi}), we obtain the same recurrence relations for the pairs $(\mu_k,\nu_k)$ and $(\tilde\mu_k,\tilde\nu_k)$
\begin{gather}
\mu_{k+2}=(2k+2)\nu_k, \qquad \nu_{k+2}=(2k+3)\mu_{k+1}+(2k+2)x\mu_k,\nonumber\\
\tilde\mu_{k+2}=(2k+2)\tilde\nu_k, \qquad \tilde\nu_{k+2}=(2k+3)\tilde\mu_{k+1}+(2k+2)x\tilde\mu_k,\label{eq:Recur_MuNu}
\end{gather}
but different initial conditions
\begin{gather*}
(\mu_0,\nu_0)=(1,0), \qquad (\mu_1,\nu_1)=(0,1),\\
(\tilde\mu_0,\tilde\nu_0)=(0,1), \qquad (\tilde\mu_1,\tilde\nu_1)=(0,0).
\end{gather*}

Now we need to reduce (\ref{eq:AsympSeries2}) to the form of (\ref{eq:AsympSeries1}). Transforming the first series of~(\ref{eq:AsympSeries2}), we have
\begin{gather*}
\sum_{k=0}^\infty \frac{\mu_k p+\nu_k}{(p^2-x)^{k+1}} =\sum_{k=0}^\infty \frac{\mu_k p+\nu_k}{p^{2k+2}}\cdot\frac{1}{(1-x/p^2)^{k+1}}\\
\hphantom{\sum_{k=0}^\infty \frac{\mu_k p+\nu_k}{(p^2-x)^{k+1}}}{}=\sum_{k=0}^\infty \frac{\mu_k p+\nu_k}{p^{2k+2}}
\sum_{m=0}^\infty \binom{m+k}{k}\frac{x^m}{p^{2m}} =\sum_{n=0}^\infty \frac{1}{p^{2n+2}}\sum_{k=0}^n \binom{n}{k}(\mu_k p+\nu_k)x^{n-k}.
\end{gather*}
For the second series of (\ref{eq:AsympSeries2}), transformation is the same. As result, we obtain expansions of~$P_n(x)$ and $Q_n(x)$ depending on parity of~$n$
\begin{gather*}
P_{2n}(x)=\sum_{k=0}^n \binom{n}{k}\mu_k(x) x^{n-k},\qquad
P_{2n+1}(x)=\sum_{k=0}^n \binom{n}{k}\nu_k(x) x^{n-k},\\
Q_{2n}(x)=\sum_{k=0}^n \binom{n}{k}\tilde\mu_k(x) x^{n-k},\qquad
Q_{2n+1}(x)=\sum_{k=0}^n \binom{n}{k}\tilde\nu_k(x) x^{n-k},
\end{gather*}
where we add an explicit dependence of the coefficients on $x$. All the coefficients satisfy difference equations of fourth order based on~(\ref{eq:Recur_MuNu}):
\begin{alignat}{3}
& \mu_k,\ \tilde\mu_k\colon \quad && y_{k+4}=(2k+3)(2k+6)y_{k+1}+(2k+2)(2k+6)xy_k, & \label{eq:Recur_Mu}\\
& \nu_k,\ \tilde\nu_k\colon \quad && y_{k+4}=(2k+4)(2k+7)y_{k+1}+(2k+2)(2k+6)xy_k.& \label{eq:Recur_Nu}
\end{alignat}
In particular,
\begin{alignat*}{3}
& \mu_k =\{1,0,0,4,12x,0,280,\ldots\}, \qquad&& \tilde\mu_k =\{0,0,2,0,0,80,120x,\ldots\},& \\
& \nu_k =\{0,1,2x,0,28,140x,120x^2,\ldots\}, \qquad && \tilde\nu_k =\{1,0,0,10,12x,0,880,\ldots\}.&
\end{alignat*}

Unfortunately, two other pairs of independent solutions of (\ref{eq:Recur_Mu}) and (\ref{eq:Recur_Nu}) are yet remain unknown, while a relation of one of the pairs (say, $\dbltilde\mu_k$ and $\dbltilde\nu_k$) to polynomials~(\ref{eq:Pol_Z}) seems quite predictable.

\section[Higher derivatives of $\Ai^2(x)$, $\Ai(x)\Bi(x)$ and $\Bi^2(x)$]{Higher derivatives of $\boldsymbol{\Ai^2(x)}$, $\boldsymbol{\Ai(x)\Bi(x)}$ and $\boldsymbol{\Bi^2(x)}$}\label{sec:Ai2_Deriv}

In a similar fashion to the evaluation of higher derivatives of Airy functions considered in Section~\ref{sec:Ai_Deriv}, let us investigate the same problem for the products of these functions. It is quite evident that the problem is reduced to a determination of polynomials $R_n(x)$, $S_n(x)$ and $T_n(x)$ satisfying the following three equations
\begin{gather}
\left(\!\!\begin{matrix} [\Ai^2(x)]^{(n)}\\ [\Ai(x)\Bi(x)]^{(n)}\\ [\Bi^2(x)]^{(n)} \end{matrix}\!\! \right)
 = \left(\!\!\begin{matrix}
\Ai^2(x) & 2\Ai(x)\Ai'(x) & \Ai'^2(x)\\
\Ai(x)\Bi(x)\!\! & \!\!\Ai(x)\Bi'(x)+\Ai'(x)\Bi(x)\!\! & \!\!\Ai'(x)\Bi'(x)\\
\Bi^2(x) & 2\Bi(x)\Bi'(x) & \Bi'^2(x)\end{matrix}\!\! \right)
\left(\! \begin{matrix} R_n(x)\\ S_n(x)\\ T_n(x) \end{matrix}\! \right)\! ,\label{eq:DerivAi^2}
\end{gather}
where $R_0(x)=1$, $S_0(x)=0$ and $T_0(x)=0$.

By taking the derivative of any equation above, we get the system of differential difference relations
\begin{gather}
R_{n+1}(x)=R'_n(x)+2xS_n(x),\nonumber\\
S_{n+1}(x)=R_n(x)+S'_n(x)+xT_n(x),\nonumber\\
T_{n+1}(x)=2S_n(x)+T'_n(x),\label{eq:DiffRecur_RST}
\end{gather}
which help to determine $R_n(x)$, $S_n(x)$ and $T_n(x)$ for any $n$ (see the Table~\ref{table2}).

\begin{table}[t]\centering
\caption{First 13 polynomials $R_n(x)$, $S_n(x)$ and $T_n(x)$.}\label{table2}\vspace{1mm}

\begin{tabular}{rccc}
 $n$\qquad & $R_n(x)$ & $S_n(x)$ & $T_n(x)$\\
 \hline
 $0$\qquad & $1$ & $0$ & $0$\medskip\\
 $1$\qquad & $0$ & $1$ & $0$\\
 $2$\qquad & $2x$ & $0$ & $2$\\
 $3$\qquad & $2$ & $4x$ & $0$\\
 $4$\qquad & $8x^2$ & $6$ & $8x$\\
 $5$\qquad & $28x$ & $16x^2$ & $20$\medskip\\
 $6$\qquad & $32x^3+28$ & $80x$ & $32x^2$\\
 $7$\qquad & $256x^2$ & $64x^3+108$ & $224x$\\
 $8$\qquad & $128x^4+728x$ & $672x^2$ & $128x^3+440$\\
 $9$\qquad & $1856x^3+728$ & $256x^4+2512x$ & $1728x^2$\\
$10$\qquad & $512x^5+10\,592x^2$ & $4608x^3+3240$ & $512x^4+8480x$\medskip\\
$11$\qquad & $11\,776x^4+27\,664x$ & $1024x^5+32\,896x^2$ &
 $11\,264x^3+14\,960$\\
$12$\qquad & $2048x^6+112\,896x^3+27\,664$ & $28\,160x^4+108\,416x$ &
 $2048x^5+99\,584x^2$
\end{tabular}
\end{table}

As before, it is better to use Airy atoms for finding the polynomials than Airy functions. By rewriting the system (\ref{eq:DerivAi^2}) as follows
\begin{gather*}
\begin{pmatrix}
\big[f^2(x)\big]^{(n)}\\ [f(x)g(x)]^{(n)}\\ \big[g^2(x)\big]^{(n)}
\end{pmatrix}
 = \begin{pmatrix}
f^2(x) & 2f(x)f'(x) & f'^2(x)\\
f(x)g(x) & f(x)g'(x)+f'(x)g(x) & f'(x)g'(x)\\
g^2(x) & 2g(x)g'(x) & g'^2(x)
\end{pmatrix}
\begin{pmatrix} R_n(x)\\ S_n(x)\\ T_n(x) \end{pmatrix}
 \end{gather*}
and applying the Wronskian $W[f,g]$ to the matrix inversion, it is easy to obtain the solution
\begin{gather}
\begin{pmatrix} R_n(x)\\ S_n(x)\\ T_n(x) \end{pmatrix}
=\begin{pmatrix}
g'^2(x) & -2f'(x)g'(x) & f'^2(x)\\
-g(x)g'(x) & f(x)g'(x)+f'(x)g(x) & -f(x)f'(x)\\
g^2(x) & -2f(x)g(x) & f^2(x)
\end{pmatrix}
\begin{pmatrix} \big[f^2(x)\big]^{(n)}\\ [f(x)g(x)]^{(n)}\\ \big[g^2(x)\big]^{(n)}
\end{pmatrix} .\label{eq:RST}
\end{gather}

Since the hypergeometric description of $f^2(x)$, $f(x)g(x)$ and $g^2(x)$ is known and simple~\cite[equation (8.3.2.38)]{Brychkov}
\begin{gather}
f^2(x)={}_1F_2\left(\frac16;\frac13,\frac23\,\Bigl|\,\frac{4x^3}{9}\right)=\sum_{k=0}^\infty \left(\frac16\right)_k \frac{12^k x^{3k}}{(3k)!},\nonumber\\
f(x)g(x)=x\cdot{}_1F_2\left(\frac12;\frac23,\frac43\,\Bigl|\,\frac{4x^3}{9}\right)=\sum_{k=0}^\infty \left(\frac12\right)_k\frac{12^k x^{3k+1}}{(3k+1)!},\nonumber\\
g^2(x)=x^2\cdot{}_1F_2\left(\frac56;\frac43,\frac53\,\Bigl|\,\frac{4x^3}{9}\right)=2\sum_{k=0}^\infty \left(\frac56\right)_k\frac{12^k x^{3k+2}}{(3k+2)!},\label{eq:AiryAtoms^2}
\end{gather}
the initial problem is reduced to manipulations of power series. Combining all three formulae of (\ref{eq:AiryAtoms^2}) together,
\begin{gather*}
\big\{f^2(x)\,\|\,f(x)g(x)\,\|\,g^2(x)\big\} =\delta!\sum_{3k+\delta\ge 0} \left(\frac{1+2\delta}{6}\right)_k\frac{12^k x^{3k+\delta}}{(3k+\delta)!},
\qquad\mbox{where}\quad \delta=\{0\,\|\,1\,\|\,2\},
\end{gather*}
(here and below, we use the symbol `$\|$' for separation of subparts) we write $n$-th derivatives of $f^2(x)$, $f(x)g(x)$ and $g^2(x)$ in the same manner:
\begin{gather*}
\big\{f^2(x)\,\|\,f(x)g(x)\,\|\, g^2(x)\big\}^{(n)} =\delta!\sum_{3k+\delta\ge n} \left(\frac{1+2\delta}{6}\right)_k \frac{12^k x^{3k+\delta-n}}{(3k+\delta-n)!},
\qquad \delta=\{0\,\|\,1\,\|\,2\},
\end{gather*}
and substitute all the series into (\ref{eq:RST}).

We begin with the case of $T_n(x)$ and employ the same approach as in Section~\ref{sec:Ai_Deriv}:
\begin{gather}
T_n(x) =g^2(x)\big[f^2(x)\big]^{(n)}-2f(x)g(x)[f(x)g(x)]^{(n)}+f^2(x)\big[g^2(x)\big]^{(n)} \label{eq:T_Series1}\\
\hphantom{T_n(x)}{}=2\sum_{3m+2\ge n} \frac{12^m x^{3m+2-n}}{(3m+2-n)!}
\sum_{0\le\delta\le 2} (-1)^\delta\sum_k \binom{3m+2-n}{3k+\delta}
\left(\frac{1+2\delta}{6}\right)_k\left(\frac{5-2\delta}{6}\right)_{m-k}.\nonumber
\end{gather}
To simplify the expression of $T_n(x)$, we shift $n$ by $2$:
\begin{gather*}
T_{n+2}(x)=2\sum_{3m\ge n} \frac{12^m x^{3m-n}}{(3m-n)!} \sum_{0\le\delta\le 2} (-1)^\delta\gamma_m(n,\delta),
\end{gather*}
where
\begin{gather*}
\gamma_m(n,\delta)=\sum_k \binom{3m-n}{3k+\delta} \left(\frac{1+2\delta}{6}\right)_k\left(\frac{5-2\delta}{6}\right)_{m-k}.
\end{gather*}
Applying Euler's beta integral we transform the inner series into an integral representation
\begin{gather*}
\gamma_m(n,\delta) =\frac{2}{\pi} m! \sin\left(\frac{\pi}{6}(1+2\delta)\right) \sum_{0\le j\le 2} \bar\omega^{\delta j}\int_0^\infty\frac{\big(1+\omega^j t^2\big)^{3m-n}}{\big(1+t^6\big)^{m+1}}\,\rmd t,\\
\sum_{0\le\delta\le 2} (-1)^\delta \gamma_m(n,\delta)=-\frac{6}{\pi}m!\Re\left\{\bar\omega\!\int_0^\infty\frac{\big(1+\omega t^2\big)^{3m-n}}{\big(1+t^6\big)^{m+1}}\,\rmd t\right\}.
\end{gather*}
As result,
\begin{gather}
T_{n+2}(x)=-\frac{1}{\pi}\sum_{3m\ge n} \frac{12^{m+1} m! x^{3m-n}}{(3m-n)!} \Re\left\{\bar\omega\int_0^\infty \frac{\big(1+\omega t^2\big)^{3m-n}}{\big(1+t^6\big)^{m+1}}\,\rmd t\right\}.\label{eq:T_Series2}
\end{gather}

Integrating the analytic function $\big(1+z^2\big)^{3m-n}/\big(1+z^6\big)^{m+1}$ over the contour $[0,R]\cup\{ R\rme^{\rmi\phi},\,\phi\in[0,\pi/3]\}\cup [R\rme^{\pi\rmi/3},0]$ and taking $R\to\infty$, one gets
\begin{gather*}
\int_0^\infty \frac{\big(1+t^2\big)^{3m-n}}{\big(1+t^6\big)^{m+1}}\,\rmd t-\bar\omega \int_\infty^0 \frac{\big(1+\omega t^2\big)^{3m-n}}{\big(1+t^6\big)^{m+1}}\,\rmd t=2\pi\rmi\mathop{\rm res}_{z=\rme^{\pi\rmi/6}}
\frac{\big(1+z^2\big)^{3m-n}}{\big(1+z^6\big)^{m+1}}.
\end{gather*}
Then the integral under consideration can be written as
\begin{align*}
\bar\omega\int_0^\infty \frac{\big(1+\omega t^2\big)^{3m-n}}{\big(1+t^6\big)^{m+1}}\,\rmd t
&=2\pi\rmi\mathop{\rm res}_{z=\rme^{\pi\rmi/6}}\frac{\big(1+z^2\big)^{3m-n}}{\big(1+z^6\big)^{m+1}}-\int_0^\infty \frac{\big(1+t^2\big)^{3m-n}}{\big(1+t^6\big)^{m+1}}\,\rmd t \\
&{}=\pi\rmi\left\{\mathop{\rm res}_{z=\rme^{\pi\rmi/6}}
-\mathop{\rm res}_{z=-\rme^{-\pi\rmi/6}}-\mathop{\rm res}_{z=\rmi}\right\}\frac{\big(1+z^2\big)^{3m-n}}{\big(1+z^6\big)^{m+1}},
\end{align*}
where we used well-known trick in the last step:
\begin{align*}
\int_0^\infty \frac{\big(1+t^2\big)^{3m-n}}{\big(1+t^6\big)^{m+1}}\,\rmd t
&=\frac12\int_\Rset \frac{\big(1+t^2\big)^{3m-n}}{\big(1+t^6\big)^{m+1}}\,\rmd t \\
&{}=\pi\rmi\left\{\mathop{\rm res}_{z=\rme^{\pi\rmi/6}}
+\mathop{\rm res}_{z=\rmi}+\mathop{\rm res}_{z=-\rme^{-\pi\rmi/6}}\right\}
\frac{\big(1+z^2\big)^{3m-n}}{\big(1+z^6\big)^{m+1}}.
\end{align*}

With standard technique on residues, it is easy to see that
\begin{gather*}
\left\{\mathop{\rm res}_{z=\rme^{\pi\rmi/6}} -\mathop{\rm res}_{z=-\rme^{-\pi\rmi/6}}\right\}\frac{\big(1+z^2\big)^{3m-n}}{\big(1+z^6\big)^{m+1}}
\end{gather*}
is a real number. Thus,
\begin{gather*}
\Re\left(\pi\rmi\left\{\mathop{\rm res}_{z=\rme^{\pi\rmi/6}} -\mathop{\rm res}_{z=-\rme^{-\pi\rmi/6}}\right\}\frac{\big(1+z^2\big)^{3m-n}}{\big(1+z^6\big)^{m+1}}\right)=0
\end{gather*}
and we need to find the residue at the pole $z=\rmi$ only. The order of the pole is $(m+1)-(3m-n)=n-2m+1$, that means $\big(1+z^2\big)^{3m-n}/\big(1+z^6\big)^{m+1}$ is holomorphic at $z=\rmi$ if $n\le 2m-1$. Therefore, nonzero summands in~(\ref{eq:T_Series2}) are possible only if $n\ge 2m$. Then the residue can be written as
\begin{gather*}
\mathop{\rm res}_{z=\rmi} \frac{\big(1+z^2\big)^{3m-n}}{\big(1+z^6\big)^{m+1}} =-\rmi\cdot h_{m,n-2m},
\end{gather*}
where
\begin{gather*}
h_{m,n}=\frac{(-1)^n}{n!} \frac{\rmd^n}{\rmd t^n} \left(\frac{1}{(t+1)^{n+1}\big(t^4+t^2+1\big)^{m+1}}\right) \biggl|_{t=1},
\end{gather*}
and (\ref{eq:T_Series2}) is reduced to the following expression:
\begin{gather}
T_{n+2}(x)=\sum_{\frac{n}{3}\le m\le\frac{n}{2}} h_{m,n-2m}\frac{12^{m+1}m!x^{3m-n}}{(3m-n)!}.\label{eq:T_Series3}
\end{gather}
In particular, $T_0(x)=T_1(x)=T_3(x)\equiv 0$ because the set of possible~$m$'s values is empty.

As before, we find $h_{m,n}$ using the generating function
\begin{align*}
H(s)=\sum_{n=0}^\infty h_{m,n}s^n&=\frac{1}{2\pi\rmi} \oint_{|z-1|=\epsilon}\sum_{n=0}^\infty\frac{(-s)^n}{\big(z^2-1\big)^{n+1}}\cdot\frac{\rmd z}{\big(z^4+z^2+1\big)^{m+1}}\\
&{}=\frac{1}{2\pi\rmi}\oint_{|z-1|=\epsilon}\frac{\rmd z}{\big(z^2-1+s\big)\big(z^4+z^2+1\big)^{m+1}}\\
&{}=\mathop{\rm res}_{z=\sqrt{1-s}}\frac{1}{\big(z^2-1+s\big)\big(z^4+z^2+1\big)^{m+1}}\\
&{}=\frac{1}{2\sqrt{1-s}\big(3-3s+s^2\big)^{m+1}}.
\end{align*}
Then
\begin{align}
h_{m,n}&=\lshad s^n\rshad H(s) =\frac{1}{2\cdot 3^{m+1}}\lshad s^n\rshad\left\{(1-s)^{-m-3/2}\left(1+\frac{s^2}{3(1-s)}\right)^{-m-1}\right\}\nonumber \\
&{}=\frac{1}{2\cdot 3^{m+1}}\,\lshad s^n\rshad\left\{\sum_{k=0}^\infty \binom{m+k}{m}\left(-\frac{s^2}{3}\right)^{k}\cdot\sum_{\ell=0}^\infty \left(m+k+\frac32\right)_{\ell}\frac{s^\ell}{\ell!}\right\}\nonumber\\
&{}=\frac{1}{2\cdot 3^{m+1}}\sum_{k=0}^{\lfloor n/2\rfloor}\binom{m+k}{m}\frac{\Gamma\bigl(m+n-k+\tfrac32\bigr)}{(n-2k)!\Gamma\bigl(k+m+\tfrac32\bigr)}\left(-\frac13\right)^{k}\nonumber\\
&{}=\frac{1}{2^{2n+1}3^{m+1}}\cdot\frac{(2m+2n+1)!m!}{(2m+1)!(m+n)!n!}\cdot{}_3F_2\left(\begin{matrix} -\frac{n}{2},-\frac{n-1}{2},m+1\\
-m-n-\frac12,m+\frac32\end{matrix} \,\biggl|\,\frac43\right).\label{eq:c_mn_1}
\end{align}

In spite of a terminating form of the hypergeometric series, the expression (\ref{eq:c_mn_1}) does not look a perfect because the argument of the series is located in the exterior of the unit disk. Having in mind the results discussed in Section~\ref{sec:HG_series}, we come back to the finite sum over~$k$ and reverse the order of summation that gives
\begin{gather*}
h_{m,2n+\delta} =\frac{(-1)^n}{2\cdot 3^{m+n+1}}\binom{m+n}{m}\left(m+n+\frac32\right)_{\delta}\\
\hphantom{h_{m,2n+\delta} =}{} \times{}_3F_2\left(\begin{matrix} -n,-m-n-\frac12,
m+n+\delta+\frac32\\ -m-n,\delta+\frac12\end{matrix}\,\biggl|\,\frac34\right).
\end{gather*}
Here we replaced the second index in $h_{m,n}$ by $2n+\delta$, where $\delta=0,1$, with the aim to avoid the appearance of many $\lfloor n/2\rfloor$'s. We recall also that if a hypergeometric function contains nonpositive integers among upper and lower parameters, then it's value is defined as
\begin{gather*}
{}_pF_q\left(\begin{matrix}-n,\ldots\\ -m-n,\ldots\end{matrix}\,\Bigl|\,
z\right)=\lim_{\epsilon\to 0} \;{}_pF_q\left(\begin{matrix} -n,\ldots\\ \epsilon-m-n,\ldots\end{matrix}\,\Bigl|\,z\right)
\end{gather*}
(see \cite[equation (2.1.4) and discussion here]{BatemanI} and \cite[Section~5.5]{GKP} for details).

Returning to (\ref{eq:T_Series3}), we obtain the required formulae for $T_n(x)$ and then for $S_n(x)$ and $R_n(x)$ on the base of relations~(\ref{eq:DiffRecur_RST}). As in Section~\ref{sec:Ai_Deriv}, we rename the coefficients in order to rewrite the final expressions in a shorter form:
\begin{gather}
T_{2n+\delta+2}(x)=\sum_{\frac{2n+\delta}{3}\le m\le n} \tilde h_{m,n,\delta}\left(\frac32,0\right)\frac{n!\bigl({-}\tfrac13\bigr)^{n-m}
2^{2m+1}x^{3m-2n-\delta}}{(n-m)!(3m-2n-\delta)!}, \label{eq:T_Series}\\
S_{2n+\delta+1}(x)=\sum_{\frac{2n+\delta}{3}\le m\le n}\tilde h_{m,n,\delta}\left(\frac12,0\right)\frac{n!\bigl({-}\tfrac13\bigr)^{n-m}2^{2m}
x^{3m-2n-\delta}}{(n-m)!(3m-2n-\delta)!},\label{eq:S_Series}\\
R_{2n+\delta}(x)=\sum_{\frac{2n+\delta}{3}\le m\le n}\left\{\tilde h_{m,n,\delta}\left(-\frac12,0\right)
-[n>0]\frac{3m-2n-\delta}{2n}\tilde h_{m,n,\delta}\left(\frac12,1\right)\right\}\nonumber\\
\hphantom{R_{2n+\delta}(x)=}{}\times \frac{n!\bigl({-}\tfrac13\bigr)^{n-m}2^{2m} x^{3m-2n-\delta}}{(n-m)!(3m-2n-\delta)!},\label{eq:R_Series}
\end{gather}
where
\begin{gather}
\tilde h_{m,n,\delta}(a,b)=(n+a)_\delta\cdot {}_3F_2\left(\begin{matrix} m-n,1-a-n,n+\delta+a\\
b-n,\delta+\frac12\end{matrix}\,\biggl|\,\frac34\right) \label{eq:RST_Coeffs}
\end{gather}
and $[n>0]$ is Iverson's symbol. In general, Iverson's symbol $[A]$ is defined for a logical state\-ment~$A$ by (see \cite[Section~2.1]{GKP})
\begin{gather*}
[A]= \begin{cases} 1, & \mbox{if $A$ is true},\\
0, & \mbox{if $A$ is false}. \end{cases}
\end{gather*}
Its presence in (\ref{eq:R_Series}) means that the second term in braces should be removed if $n=0$.

Now we briefly discuss two corollaries of the above formulae~-- a difference equation of third order and special values of the function ${}_3F_2$.

A difference equation for polynomials $R_n(x)$, $S_n(x)$ and $T_n(x)$ follows by applying the way used in Section~\ref{sec:HG_series}. It is known that products of Airy atoms $f^2(x)$, $f(x)g(x)$ and $g^2(x)$ satisfy the equation $y'''-4xy'-2y=0$. Since (\ref{eq:RST}), the generating functions of all three polynomial families,
\begin{gather*}
{\bf R}(x,t)=\sum_{n=0}^\infty R_n(x)\frac{t^n}{n!},\qquad
{\bf S}(x,t)=\sum_{n=0}^\infty S_n(x)\frac{t^n}{n!},\qquad
{\bf T}(x,t)=\sum_{n=0}^\infty T_n(x)\frac{t^n}{n!},
\end{gather*}
are solutions of the differential equation
\begin{gather*}
\left\{\frac{\rmd^3}{\rmd t^3}-4(x+t)\frac{\rmd}{\rmd t}-2\right\} {\bf Y}(x,t)=0.
\end{gather*}
Expanding a solution in a power series,
\begin{gather*}
{\bf Y}(x,t)=\sum_{n=0}^\infty Y_n(x)\frac{t^n}{n!},
\end{gather*}
one can find a recurrence relation for its coefficients, functions $Y_n(x)$:
\begin{gather}
Y_{n+3}(x)=4xY_{n+1}(x)+(4n+2)Y_n(x).\label{eq:Recur_RST}
\end{gather}
Thus, in contrast with (\ref{eq:Recur_PQZ}), three independent solutions of (\ref{eq:Recur_RST}) are known:
\begin{alignat*}{5}
& R_n(x)\colon \quad && R_0(x)=1, \qquad && R_1(x)=0, \qquad && R_2(x)=2x,&\\
& S_n(x)\colon \quad && S_0(x)=0, \qquad && S_1(x)=1, \qquad && S_2(x)=0,&\\
& T_n(x)\colon \quad && T_0(x)=0, \qquad && T_1(x)=0, \qquad && T_2(x)=2,&
\end{alignat*}
and the general solution of (\ref{eq:Recur_RST}) with arbitrary initial conditions is
\begin{gather*}
Y_n(x)=Y_0(x)R_n(x)+Y_1(x)S_n(x)+\left(\frac{Y_2(x)}{2}-xY_0(x)\right)T_n(x).
\end{gather*}

Relations containing special values of the function ${}_3F_2$ are based on expansions of polynomials $R_n(x)$, $S_n(x)$ and $T_n(x)$ for small $x$'s. For the case of $T_n(x)$, we can use (\ref{eq:T_Series1}) and obtain
\begin{gather*}
T_{3n}(x)=12^n\left\{\left(\frac56\right)_n -2\left(\frac12\right)_n+\left(\frac16\right)_n\right\}x^2+\cdots,\\
T_{3n+1}(x)=2\cdot 12^n\left\{\left(\frac56\right)_{n}-\left(\frac12\right)_{n}\right\}x\\
\hphantom{T_{3n+1}(x)=}{} +\frac{12^n}{2}\left\{(2n+3)\left(\frac56\right)_{n}-(8n+5)\left(\frac12\right)_{n}+2\left(\frac76\right)_{n}\right\}x^4+\cdots,\\
T_{3n+2}(x)=2\cdot 12^n\left(\frac56\right)_n +2\cdot 12^n\left\{(2n+2)\left(\frac56\right)_n -3\left(\frac32\right)_n+\left(\frac76\right)_n\right\}x^3+\cdots.
\end{gather*}
For two other polynomial families, similar expansions can be found either by reducing to $T_n(x)$ with (\ref{eq:DiffRecur_RST}) or by applying (\ref{eq:RST}) directly. The results are as follows:
\begin{gather*}
S_{3n}(x)=12^n\left\{\left(\frac12\right)_n-\left(\frac16\right)_n\right\}x+\cdots,\\
S_{3n+1}(x)=12^n\left(\frac12\right)_n+12^n\left\{(2n+2)\left(\frac12\right)_n-\left(\frac56\right)_n-\left(\frac76\right)_n\right\}x^3+\cdots,\\
S_{3n+2}(x)=12^n\left\{3\left(\frac32\right)_n-\left(\frac56\right)_n-2\left(\frac76\right)_n\right\}x^2+\cdots,\\
R_{3n}(x)=12^n\left(\frac16\right)_n+12^n\left\{(2n+1)\left(\frac16\right)_n-\left(\frac12\right)_n\right\}x^3 +\cdots,\\
R_{3n+1}(x)=12^n\left\{\left(\frac76\right)_n-\left(\frac12\right)_n\right\}x^2+\cdots,\\
R_{3n+2}(x)=12^{n+1}\left(\frac16\right)_{n+1}x+\frac{12^n}{2}\left\{(2n+5)\left(\frac76\right)_n+\left(\frac56\right)_n-6\left(\frac32\right)_n\right\}x^4+\cdots.
\end{gather*}
Then, using (\ref{eq:RST_Coeffs}), we can write special values of the function ${}_3F_2$ in terms of the Pochhammer symbol. As in Section~\ref{sec:HG_series}, Carlson's theorem and asymptotical behaviour of gamma function
help to replace an integer $n$ by a real $a$ (see Appendix~\ref{appendixA} for details). Below we present the formulae corresponding to the shortest cases of the polynomials' coefficients only, namely, $T_{3n+2}(0)$, $S_{3n+1}(0)$, $R_{3n}(0)$ and $R'_{3n+2}(0)$. The first formula is for $\delta=0$, the second~-- for $\delta=1$.

The case $T_{3n+2}(0)=2\cdot 12^n\bigl(\tfrac56\bigr)_n$ gives
\begin{gather}
{}_3F_2\left(\begin{matrix} a,3a-\frac12,\frac32-3a\\
3a,\frac12\end{matrix}\,\Bigl|\,\frac34\right) =\frac{2\Gamma\bigl(\tfrac16\bigr)\Gamma(3a)\sin\bigl(\tfrac{\pi}{6}+\pi a\bigr)}{3^{3a}\Gamma\bigl(2a+\tfrac16\bigr)\Gamma(a)},\label{eq:HyperG_3F2_x=3/4_a}\\
{}_3F_2\left(\begin{matrix} a,3a-\frac32,\frac72-3a\\ 3a-1,\frac32\end{matrix}\, \Bigl|\,\frac34\right) =\frac{5-12a}{(1-3a)(5-6a)}\cdot\frac{2\Gamma\bigl(\tfrac16\bigr)\Gamma(3a)\sin\bigl(\tfrac{\pi}{6}+\pi a\bigr)}
{3^{3a}\Gamma\bigl(2a+\tfrac16\bigr)\Gamma(a)}.\label{eq:HyperG_3F2_x=3/4_b}
\end{gather}

The case $S_{3n+1}(0)=12^n\bigl(\tfrac12\bigr)_n$ gives
\begin{gather}
{}_3F_2\left(\begin{matrix} a,\frac12-3a,\frac12+3a\\ 3a,\tfrac12\end{matrix}\, \Bigl|\,\frac34\right) =\frac{4\sqrt\pi\Gamma(3a)\cos\pi a}{3^{3a}\Gamma\bigl(2a+\tfrac12\bigr)\Gamma(a)},\label{eq:HyperG_3F2_x=3/4_c}\\
{}_3F_2\left(\begin{matrix} a,3a-\frac12,\frac52-3a\\ 3a-1,\tfrac32\end{matrix}\, \Bigl|\,\frac34\right) =\frac{1-4a}{(1-2a)(1-3a)}\cdot
\frac{4\sqrt\pi\Gamma(3a)\cos\pi a} {3^{3a}\Gamma\bigl(2a+\tfrac12\bigr)\Gamma(a)}.\label{eq:HyperG_3F2_x=3/4_d}
\end{gather}

The case $R_{3n}(0)=12^n\bigl(\tfrac16\bigr)_n$ gives
\begin{gather*}
{}_3F_2\left(\begin{matrix} a,3a+\frac32,-\frac12-3a\\ 3a,\tfrac12\end{matrix}\, \Bigl|\,\frac34\right) =\frac{2\Gamma\bigl(\tfrac56\bigr)\Gamma(3a) \sin\bigl(\tfrac{\pi}{6}-\pi a\bigr)}
{3^{3a}\Gamma\bigl(2a+\tfrac56\bigr)\Gamma(a)},\\ 
{}_3F_2\left(\begin{matrix} a,3a+\frac12,\frac32-3a\\ 3a-1,\tfrac32\end{matrix}\,\Bigl|\,\frac34\right) =\frac{1-12a}{(1-3a)(1-6a)}\cdot\frac{2\Gamma\bigl(\tfrac56\bigr)\Gamma(3a)
\sin\bigl(\tfrac{\pi}{6}-\pi a\bigr)}{3^{3a}\Gamma\bigl(2a+\tfrac56\bigr)\Gamma(a)}.
\end{gather*}

The case $R'_{3n+2}(0)=12^{n+1}\bigl(\tfrac16\bigr)_{n+1}$ gives
\begin{gather}
{}_3F_2\left(\begin{matrix} a,3a+\frac12,\frac12-3a\\ 3a-1,\frac12\end{matrix}\, \Bigl|\,\frac34\right)\nonumber\\
\qquad{} {}=\frac{\Gamma(3a)}{(1-3a)3^{3a}\Gamma(a)}\left\{ \frac{\Gamma\bigl(\tfrac16\bigr)\sin\bigl(\tfrac{\pi}{6}+\pi a\bigr)} {\Gamma\bigl(2a+\tfrac16\bigr)}+(1-12a)
\frac{\Gamma\bigl(\tfrac56\bigr)\sin\bigl(\tfrac{\pi}{6}-\pi a\bigr)} {\Gamma\bigl(2a+\tfrac56\bigr)}\right\},\nonumber \\ 
{}_3F_2\left(\begin{matrix} a,3a-\frac12,\frac52-3a\\ 3a-2,\frac32\end{matrix}\, \Bigl|\,\frac34\right) =\frac{\Gamma(3a)}{(1-2a)(1-3a)(2-3a)3^{3a+1}\Gamma(a)}\nonumber\\
\qquad {}\times\left\{(5-12a)\frac{\Gamma\bigl(\tfrac16\bigr) \sin\bigl(\tfrac{\pi}{6}+\pi a\bigr)}{\Gamma\bigl(2a+\tfrac16\bigr)}
+(1-12a)(7-12a)\frac{\Gamma\bigl(\tfrac56\bigr)\sin\bigl(\tfrac{\pi}{6}-\pi a\bigr)}{\Gamma\bigl(2a+\tfrac56\bigr)}\right\}.\label{eq:HyperG_3F2_x=3/4_h}
\end{gather}

Finally, it is worth noting that $R_n(x)$, $S_n(x)$ and $T_n(x)$ can be expressed in terms of $P_n(x)$ and $Q_n(x)$. Since
\begin{gather*}
\bigl[\Ai^2(x)\bigr]^{(n)} =\sum_{k=0}^n \binom{n}{k}\Ai^{(k)}(x)\Ai^{(n-k)}(x),
\end{gather*}
we have that
\begin{gather*}
\bigl\{R_n(x)\,\big\|\, 2S_n(x)\, \big\|\, T_n(x)\bigr\}\\
\qquad {}=\sum_{k=0}^n \binom{n}{k} \bigl\{ P_k(x)P_{n-k}(x)\,\big\|\,
P_k(x)Q_{n-k}(x)+Q_k(x)P_{n-k}(x)\,\big\|\, Q_k(x)Q_{n-k}(x)\bigr\}.
\end{gather*}
However, it seems difficult to obtain the expansions (\ref{eq:T_Series})--(\ref{eq:R_Series}) by using formulae (\ref{eq:Solution_Q}) and~(\ref{eq:Solution_P}) directly.

\section{Concluding remarks}\label{sec:Extro}

In this paper, we have found higher derivatives of Airy functions and their products in a closed form. Explicit formulae for the polynomials which are contained in these derivatives help to obtain special values of hypergeometric functions ${}_2F_1$ and ${}_3F_2$. We did not try to give a complete solution of this problem and considered only a few simple cases. In our view, a construction of the general solution will not require significant efforts. Moreover, combining this approach with the ideas presented in~\cite{Gessel}, one can find more general values of ${}_3F_2$-function depending on two parameters. For example (cf.\ with (\ref{eq:HyperG_3F2_x=3/4_c}) and (\ref{eq:HyperG_3F2_x=3/4_d})),
\begin{gather*}
{}_3F_2\left(\begin{matrix} b,\frac12-3a,\frac12+3a\\
3b,\frac12\end{matrix}\,\Bigl|\,\frac34\right) =\frac{4\Gamma\bigl(\tfrac12+a-b\bigr)\Gamma(3b)}
{3^{3b}\Gamma\bigl(\tfrac12+a+b\bigr)\Gamma(b)}\cos\pi a\cos[\pi(b-a)],\\ 
{}_3F_2\left(\begin{matrix} b,1-3a,1+3a\\ 3b-1,\frac32\end{matrix}\, \Bigl|\,\frac34\right)
=\frac{4\Gamma(1+a-b)\Gamma(3b-1)} {3^{3b}a\Gamma(a+b)\Gamma(b)}\sin\pi a \sin[\pi(b-a)].
\end{gather*}

Nevertheless, there are some problems which are yet remain open and look quite difficult. For example, how to find the solution $Z_n(x)$ of the equation~(\ref{eq:Recur_PQZ}) and to reveal its relation to Airy functions.

The next problem is connected with zeros of all above polynomials. Let us define {\it reduced polynomials} removing a zero at the origin and an excessive cubicity:
\begin{alignat*}{3}
& P_{3n+\delta}(x)=x^{\{0\,\|\,2\,\|\,1\}}\tilde P_{3n+\delta}\big(x^3\big), \qquad && R_{3n+\delta}(x)=x^{\{0\,\|\,2\,\|\,1\}}\tilde R_{3n+\delta}\big(x^3\big),&\\
& Q_{3n+\delta}(x)=x^{\{1\,\|\,0\,\|\,2\}}\tilde Q_{3n+\delta}\big(x^3\big), \qquad && S_{3n+\delta}(x)=x^{\{1\,\|\,0\,\|\,2\}}\tilde S_{3n+\delta}\big(x^3\big),&\\
& Z_{3n+\delta}(x)=x^{\{2\,\|\,1\,\|\,0\}}\tilde Z_{3n+\delta}\big(x^3\big), \qquad && T_{3n+\delta}(x)=x^{\{2\,\|\,1\,\|\,0\}}\tilde T_{3n+\delta}\big(x^3\big),&
\end{alignat*}
where $\delta=\{0\,\|\,1\,\|\,2\}$. For example,
\begin{gather*}
\tilde Q_{15}(x)=x^2+770x+8680,\qquad \tilde R_{12}(x)=2048x^2+112\,896x+27\,664.
\end{gather*}
Numerical investigations show that zeros of all reduced polynomials are real (therefore, negative) and simple but we could not prove this and find an asymptotical behaviour of zeros and extreme points of the polynomials.

In conclusion, it is interesting to note that negative counterparts of the coefficients (\ref{eq:Q_Coeffs}) appear in some integral transforms of Airy functions. For example,
\begin{gather*}
\int_0^\infty x^{3n+1/2}\Ai^2(ax)\{\Bi(bx)-\rmi\,\Ai(bx)\}\,\rmd x=\frac{2^{1/3}}{2^{2n+2}3\pi n!a^{3n+3/2}}\\
{}\times\sum_{m=0}^{2n}\! \left(\frac12\right)_{m}\! (3n-m)!\lshad t^{2n-m}\rshad\!\left(1-t+\frac{t^2}{3}\right)^{n}\!
\left\{\Lambda_{m,3n-m}\left(\frac{b}{2^{2/3}a}\right)\!-\omega\Lambda_{m,3n-m}\left(\frac{\omega b}{2^{2/3}a}\right)\!\right\},
\end{gather*}
where $a>0$, $b<2^{2/3}a$, $\Re(\sqrt{2^{2/3}a-\omega b})>0$, and
\begin{gather*}
\Lambda_{n,N}(t)=\frac{1}{t^{N+1}}\left\{\frac{1}{(1-t)^{n+1/2}}
-\sum_{k=0}^N \left(n+\frac12\right)_{k}\frac{t^k}{k!}\right\}
=\frac{1}{t^{N+1}}\sum_{k=N+1}^\infty \left(n+\frac12\right)_{k}\frac{t^k}{k!}.
\end{gather*}

\appendix
\section[Appendix: Evaluation of the ${}_3F_2$-function special value]{Appendix: Evaluation of the $\boldsymbol{{}_3F_2}$-function special value}\label{appendixA}

Here, we prove the relation (\ref{eq:HyperG_3F2_x=3/4_a}).

Since $T_{3n+2}(0)=2\cdot 12^n\bigl(\tfrac56\bigr)_n$, it follows that $T_{6n+3\delta+2}(0)=2\cdot 12^{2n+\delta}\bigl(\tfrac56\bigr)_{2n+\delta}$. We can also obtain $T_{6n+3\delta+2}(0)=T_{2[3n+\delta]+\delta+2}(0)$ from the expansion~\eqref{eq:T_Series}. Equating both expressions leads to
\begin{gather*}
2\cdot 12^{2n+\delta}\left(\frac56\right)_{2n+\delta}=\tilde h_{2n+\delta,3n+\delta,\delta}\left(\frac32,0\right) \frac{2^{4n+2\delta+1}(3n+\delta)!}{(-3)^n n!}.
\end{gather*}
Then
\begin{gather*}
\left(3n+\delta+\frac32\right)_{\delta}\cdot{}_3F_2\left(\begin{matrix}-n,-\frac12-\delta-3n,3n+2\delta+\frac32 \\ -\delta-3n,\delta+\frac12\end{matrix}\,\Bigl|\,\frac34\right)\\
\qquad{}=\frac{n!}{(3n+\delta)!}\sum_{k=0}^n \frac{\Gamma\bigl(3n+2\delta+\tfrac32+k\bigr)(3n+\delta-k)!(-3)^k}{\Gamma\bigl(3n+\delta+\tfrac32-k\bigr)(n-k)!(2k+\delta)!}
=\frac{(-1)^n3^{3n+\delta}\bigl(\tfrac56\bigr)_{2n+\delta}n!}{(3n+\delta)!}.
\end{gather*}
Below we restrict ourselves to the case $\delta=0$ only:
\begin{gather*}
{}_3F_2\left(\begin{matrix} -n,-\frac12-3n,\frac32+3n \\ -3n,\frac12\end{matrix}\,\Bigl|\,\frac34\right) =\frac{(-1)^n3^{3n}\bigl(\tfrac56\bigr)_{2n}n!}{(3n)!}.
\end{gather*}

As mentioned above, the case when a hypergeometric function contains nonpositive integers among its upper and lower parameters is special and cannot be generalized to deal with real values. Nevertheless, let us replace $n$ by $-a$ and consider a hypothetical identity
\begin{gather*}
{\bf F}(a) ={}_3F_2\left(\begin{matrix} a,3a-\frac12,\frac32-3a \\ 3a,\frac12\end{matrix}\,\Bigl|\,\frac34\right)
=\tau(a)\cdot\frac{\Gamma\bigl(\tfrac56-2a\bigr)\Gamma(1-a)} {3^{3a}\Gamma\bigl(\tfrac56\bigr)\Gamma(1-3a)},
\end{gather*}
where ${\bf F}(a)$ is a short-hand notation for our hypergeometric function, and $\tau(a)$ is an unknown function.

\begin{figure}[t]\centering
\makebox[0.5cm]{(a)~}
\includegraphics[width=0.9\textwidth,keepaspectratio]{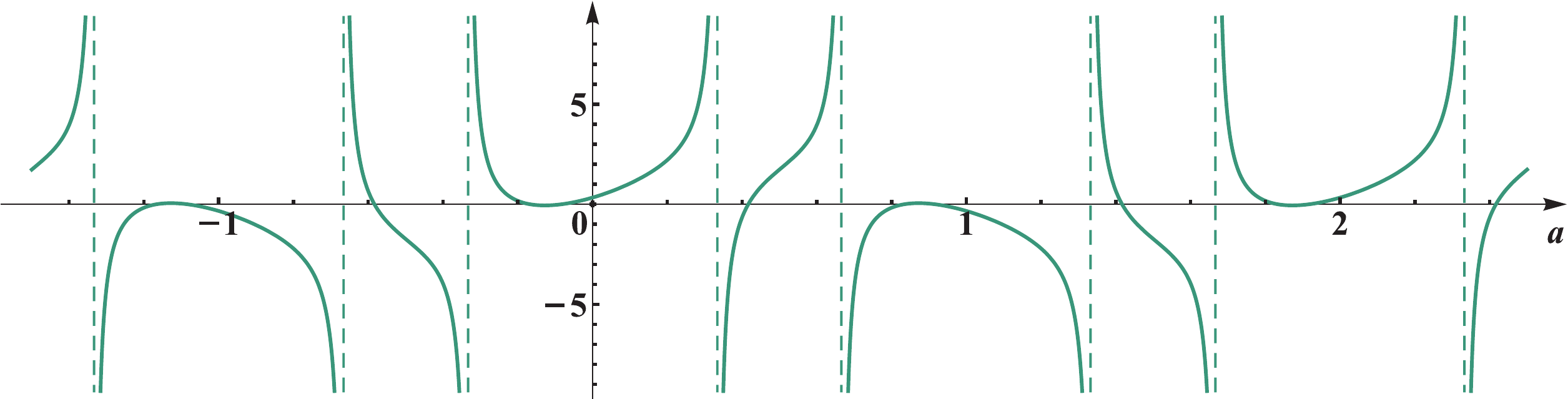}

\vspace{2pt}

\makebox[0.5cm]{(b)~}
\includegraphics[width=0.9\textwidth,keepaspectratio]{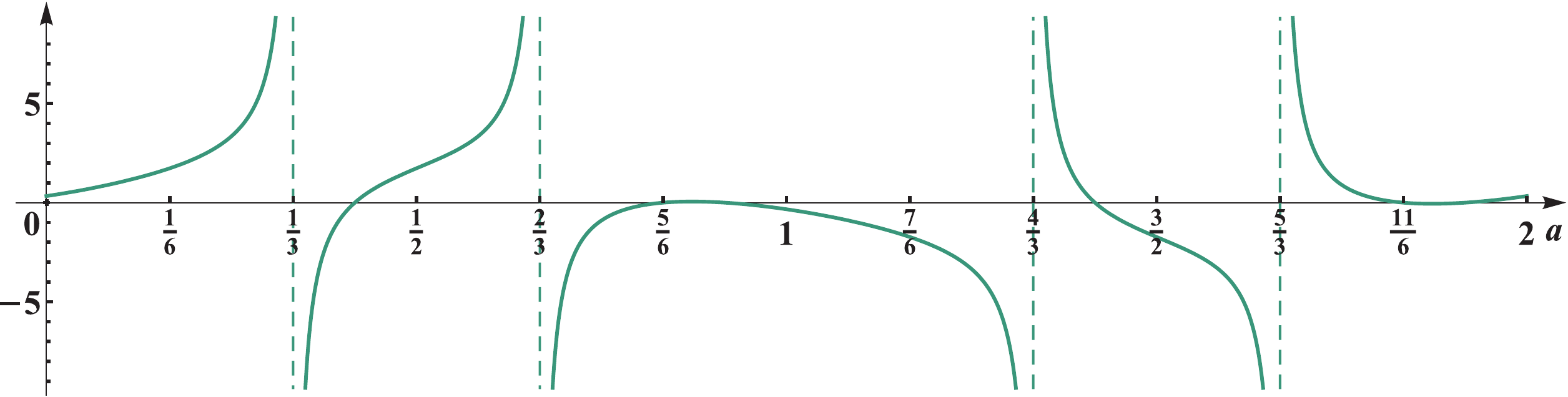}

\vspace{2pt}

\makebox[0.5cm]{(c)~}%
\includegraphics[width=0.9\textwidth,keepaspectratio]{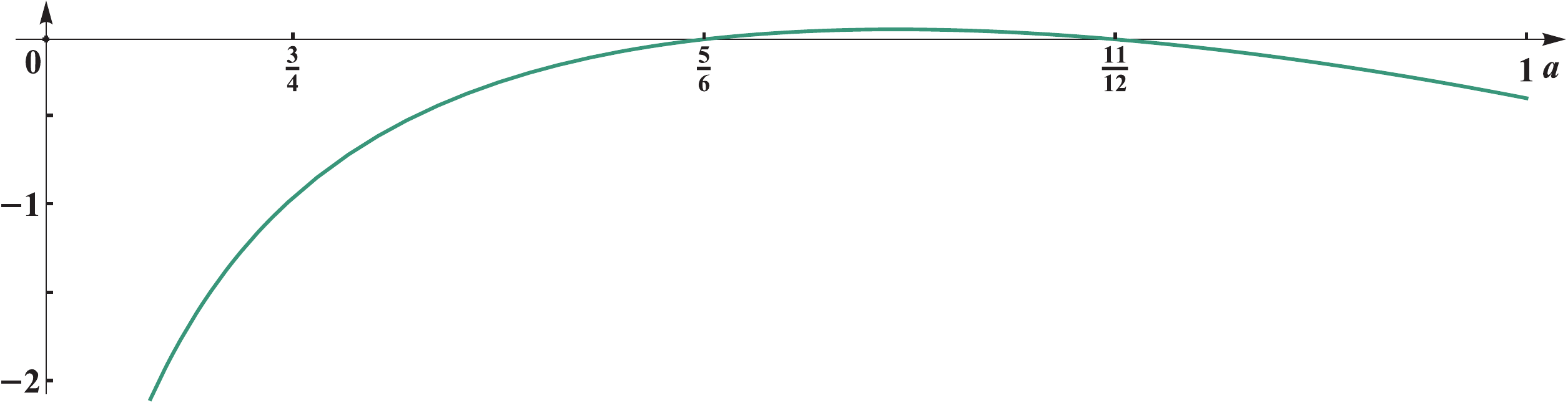}

\caption{The function $\tau(a)$ for $a\in[-1.5,2.5]$~(a), $a\in[0,2]$~(b), and $a\in[0.7,1]$~(c).}\label{fig:TauFunc}
\end{figure}

The numerical plotting procedure shows that the function
\begin{gather*}
\tau(a)={\bf F}(a)\cdot\frac{3^{3a}\Gamma\bigl(\tfrac56\bigr)\Gamma(1-3a)} {\Gamma\bigl(\tfrac56-2a\bigr)\Gamma(1-a)}
\end{gather*}
is a periodic function with period~2 (see Fig.~\ref{fig:TauFunc}). A closer scrutiny of the curve reveals that
\begin{gather*}
\tau\left(n+\frac13\right)=\tau\left(n+\frac23\right)=\infty,\qquad \tau\left(n+\frac{5}{12}\right)=\tau\left(n+\frac56\right)
=\tau\left(n+\frac{11}{12}\right)=0
\end{gather*}
for all $n\in\Zset$.

\begin{figure}[t]\centering
\includegraphics[width=0.9\textwidth,keepaspectratio]{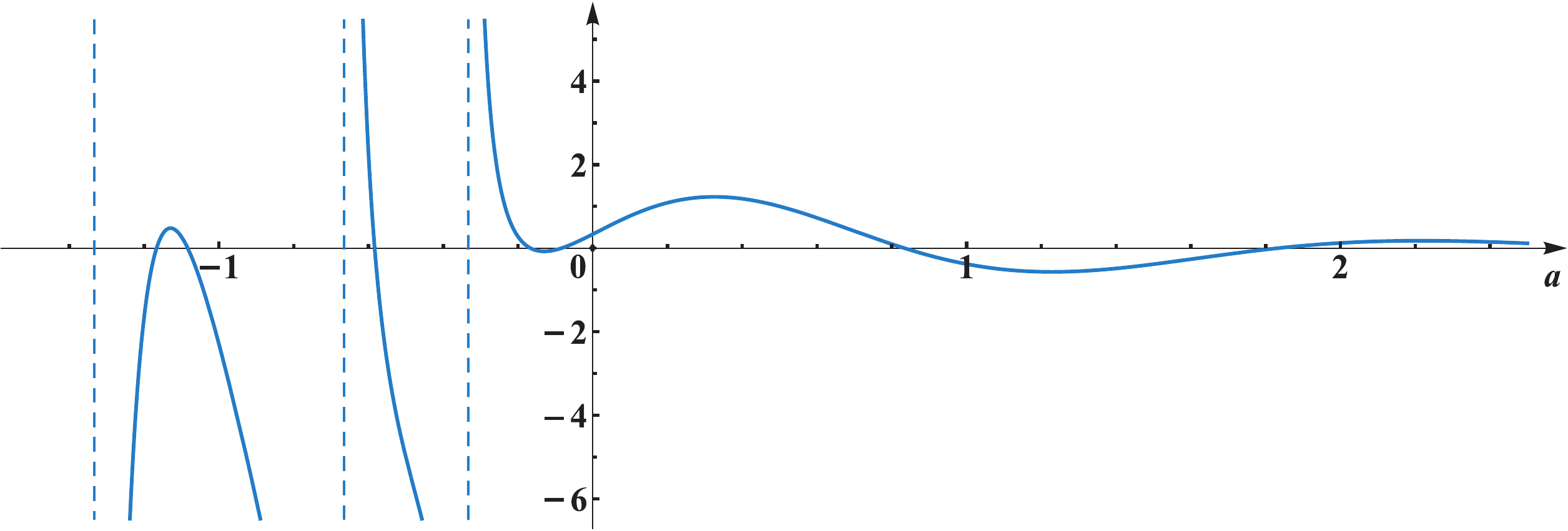}
\caption{The hypergeometric function ${\bf F}(a)$ for $a\in[-1.5,2.5]$.}\label{fig:3F2}
\end{figure}

Combining all the above zeros and poles together, let us introduce the function
\begin{gather*}
\tilde\tau(a)=\frac{\sin\bigl(\pi\bigl[a-\tfrac5{12}\bigr]\bigr) \sin\bigl(\pi\bigl[a-\tfrac56\bigr]\bigr)
\sin\bigl(\pi\bigl[a-\tfrac{11}{12}\bigr]\bigr)}{\sin\bigl(\pi\bigl[a-\tfrac13\bigr]\bigr)
\sin\bigl(\pi\bigl[a-\tfrac23\bigr]\bigr)}=-\frac{\sin\bigl(\pi\bigl[a-\tfrac56\bigr]\bigr)
\sin\bigl(\pi\bigl[2a-\tfrac56\bigr]\bigr)}{2\sin\bigl(\pi\bigl[a-\tfrac13\bigr]\bigr)\sin\bigl(\pi\bigl[a-\tfrac23\bigr]\bigr)}
\end{gather*}
with a hope that $\tilde\tau(a)=\tau(a)$. However, the numerical plotting procedure reveals that
\begin{gather}
\frac{{\bf F}(a)} {\tilde\tau(a)\cdot\dfrac{\Gamma\bigl(\tfrac56-2a\bigr)\Gamma(1-a)}{3^{3a}\Gamma\bigl(\tfrac56\bigr)\Gamma(1-3a)}}=-2\label{eq:Constant}
\end{gather}
for all $a$. As result, one can obtain the formula
\begin{gather*}
{\bf F}(a) =\frac{\sin\bigl(\pi\bigl[a-\tfrac56\bigr]\bigr)\sin\bigl(\pi\bigl[2a-\tfrac56\bigr]\bigr)}{\sin\bigl(\pi\bigl[a-\tfrac13\bigr]\bigr)
\sin\bigl(\pi\bigl[a-\tfrac23\bigr]\bigr)}\cdot\frac{\Gamma\bigl(\tfrac56-2a\bigr)\Gamma(1-a)}{3^{3a}\Gamma\bigl(\tfrac56\bigr)\Gamma(1-3a)}\\
 \hphantom{{\bf F}(a)}{} =\frac{\Gamma\bigl(\tfrac16\bigr)\Gamma\bigl(a+\tfrac13\bigr)
\Gamma\bigl(a+\tfrac23\bigr)}{\pi\sqrt{3}\,\Gamma\bigl(2a+\tfrac16\bigr)}\cdot\sin\Bigl(\pi a+\frac{\pi}{6}\Bigr),
\end{gather*}
which coincides with equation~(54) (see Fig.~\ref{fig:3F2}, where ${\bf F}(a)$ is shown for $a\in[-1.5,2.5]$).

Now, ``when you have satisfied yourself that the theorem is true, you start proving it'' \cite[Chapter~5]{PolyaI}.
At this stage, it is convenient to replace the parameter $a$ by a complex $z$:
\begin{gather*}
{\bf F}(z)={}_3F_2\left(\begin{matrix} z,3z-\frac12,\frac32-3z \\ 3z,\frac12\end{matrix}\,\Bigl|\,\frac34\right)\stackrel{?}{=}{\bf F}_0(z)
=\frac{\Gamma\bigl(\tfrac16\bigr)\Gamma\bigl(z+\tfrac13\bigr) \Gamma\bigl(z+\tfrac23\bigr)}{\pi\sqrt{3}\Gamma\bigl(2z+\tfrac16\bigr)}\sin\left(\pi z+\frac{\pi}{6}\right).
\end{gather*}
We will prove that ${\bf F}(z)={\bf F}_0(z)$ for all $z$ in the region of joint analyticity by fitting ${\bf F}(z)-{\bf F}_0(z)$ to conditions of Carlson's theorem.

First of all, we cannot use the left half-plane, $\Re z\le 0$, since ${\bf F}(z)$ is meromorphic there.

Let us consider the right half-plane, $\Re z\ge 0$, where, as we will see later, ${\bf F}(z)$ is holomorphic. We need to define a sequence of positive numbers $\{z_n\}$ with $z_n\to\infty$ such that ${\bf F}(z_n)$ can be found easily in order to examine the validity of the relation ${\bf F}(z_n) ={\bf F}_0(z_n)$. The simplest way is to choose $z_n$'s for which the hypergeometric series {\it terminates}: $z=\tfrac16$ and $\tfrac32-3z=0,-1,-2,\ldots$. Thus, one can obtain the following sequences:
\begin{gather*}
z_n=n+\frac16,\qquad \tilde z_n=n+\frac12,\qquad \dbltilde z_n=n+\frac56,\qquad n\ge 0.
\end{gather*}

Since ${\bf F}_0(\dbltilde z_n)=0$, we consider the sequence $\{\dbltilde z_n\}$ first. Then
\begin{gather*}
{\bf F}(\dbltilde z_n) ={}_3F_2\left(\begin{matrix} n+\frac56,3n+2,-1-3n \\ 3n+\frac52,\frac12\end{matrix}\,\Bigl|\,\frac34\right)
=\sum_{k=0}^{3n+1} \binom{3n+1+k}{2k} \frac{\bigl(n+\tfrac56\bigr)_k}{\bigl(3n+\tfrac52\bigr)_k} (-3)^k.
\end{gather*}
We use Zeilberger's algorithm (see \cite[Chapter~5]{GKP} and \cite{Zeilberger} for detailed description) which for a~given finite sum of hypergeometric terms, $F_n=\sum_k f(n,k)$, tries to construct a linear difference operator annihilating $F_n$, ${\cal L}(n)F_n=0$, and a rational function $R(n,k)$ such that
\begin{gather*}
{\cal L}(n)f(n,k)=G(n,k+1)-G(n,k),\qquad\mbox{where}\quad G(n,k)=R(n,k)f(n,k).
\end{gather*}
There are various implementations of the algorithm in modern computer algebra systems, for example, the Maple package {\sf EKHAD}. Applying it to ${\bf F}(\dbltilde z_n)$, we have found the following pair:
\begin{gather*}
{\cal L}(n)={\cal L}(n\,|\,\dbltilde z_n) =(12n+11)(12n+17){\cal N}+(6n+7)(6n+9),\\
R(n,k)=R(n,k\,|\,\dbltilde z_n)=12k(2k-1)\big(12n^2+32n+21\big)\\
\hphantom{R(n,k)=}{}\times\frac{2k^3-18k^2(n+1)-2k\big(81n^2+153n+73\big)-3(n+1)\big(90n^2+162n+71\big)}
{(6n+7+2k)(6n+5+2k)(3n+2-k)(3n+3-k)(3n+4-k)},
\end{gather*}
where $\cal N$ denotes the forward shift operator, i.e., ${\cal N}F_n=F_{n+1}$. Since ${\cal L}(n){\bf F}_0(\dbltilde z_n)=0$ and ${\bf F}(\dbltilde z_0)=0$, we immediately obtain that ${\bf F}(\dbltilde z_n)={\bf F}_0(\dbltilde z_n)$ for all~$n$.

Using this approach, one can find $({\cal L},R)$-pairs for two other sequences:
\begin{gather*}
{\bf F}(\tilde z_n)
={}_3F_2\left(\begin{matrix} n+\frac12,3n+1,-3n \\
3n+\frac32,\frac12\end{matrix}\,\Bigl|\,\frac34\right)\qquad\mbox{and}\qquad
{\bf F}(z_n)={}_3F_2 \left(\begin{matrix} n+\frac16,3n,1-3n \\
3n+\frac12,\frac12\end{matrix}\,\Bigl|\,\frac34\right).
\end{gather*}
Since $\tilde z_n=\dbltilde z_n-\frac13$ and $z_n=\dbltilde z_n-\frac23$, it is not surprising that
\begin{alignat*}{3}
& {\cal L}(n\,|\,\tilde z_n) ={\cal L}\left(n-\frac13\,\Bigl|\,\dbltilde z_n\right),\qquad && R(n,k\,|\,\tilde z_n)=R\left(n-\frac13,k\,\Bigl|\,\dbltilde z_n\right),& \\
& {\cal L}(n\,|\,z_n) ={\cal L}\left(n-\frac23\,\Bigl|\,\dbltilde z_n\right),\qquad && R(n,k\,|\,z_n)=R\left(n-\frac23,k\,\Bigl|\,\dbltilde z_n\right).&
\end{alignat*}

The final part of the proof of the conjectured relations, ${\bf F}(\tilde z_n)={\bf F}_0(\tilde z_n)$ and ${\bf F}(z_n)={\bf F}_0(z_n)$, is completely routine, by checking that
\begin{gather*}
{\bf F}_0(\tilde z_n)=(-1)^n \frac{\bigl(\tfrac56\bigr)_n\bigl(\tfrac76\bigr)_n} {\bigl(\tfrac76\bigr)_{2n}}\qquad\mbox{and}\qquad
{\bf F}_0(z_n)=(-1)^n \frac{\bigl(\tfrac12\bigr)_n\bigl(\tfrac56\bigr)_n} {\bigl(\tfrac12\bigr)_{2n}}
\end{gather*}
vanish under the action of the corresponding ${\cal L}$-operator, and by checking the trivial initial condition for $n=0$.

Now, we turn to Carlson's theorem. The theorem is a special case of another theorem of Carlson (see \cite[Section~5.8]{Titchmarsh_TheorFunc_ENG}): {\it Let $f(z)$ be regular and of the form ${\cal O}\big(\rme^{c|z|}\big)$ for $\Re z\ge 0$; and let $f(z)={\cal O}\big(\rme^{-\alpha|z|}\big)$, where $\alpha>0$, on the imaginary axis. Then $f(z)=0$ identically.} We will refer to it as Theorem~5.8. In \cite{Titchmarsh_TheorFunc_ENG}, Carlson's theorem has been proved in three steps by checking that $F(z)=f(z)/\sin\pi z$ satisfies the conditions of Theorem~5.8. First, the function $1/\sin\pi z$ is bounded on the circles $|z|=n+\tfrac12$. Hence $F(z)={\cal O}\big(\rme^{c|z|}\big)$ on these circles, and also on the imaginary axis. Second, since $F(z)$ is regular, it follows that for $z\in\bigl\{\Re z\ge 0,\, n-\tfrac12\le|z|\le n+\tfrac12\bigr\}$
\begin{gather*}
F(z)={\cal O}\big(\rme^{c(n+1/2)}\big)={\cal O}\big(\rme^{c|z|}\big)
\end{gather*}
and so $F(z)$ is of the form throughtout the right half-plane. Third,
\begin{gather*}
\bigl|F(\rmi y)\bigr|=\frac{\bigl|f(\rmi y)\bigr|}{\sinh \pi|y|} ={\cal O}\big(\rme^{-(\pi-c)|y|}\big),
\end{gather*}
and, therefore, the identity $F(z)=0$ follows from Theorem~5.8.

The next three statements are simple corollaries of Theorem~5.8 and can be proved in the same way as Carlson's theorem itself.

\begin{Corollary}\label{Corollary1} If $f(z)$ is regular and of the form ${\cal O}\big(\rme^{c|z|}\big)$, where $c<2\pi$, for $\Re z\ge 0$, and $f(z)=0$ for $z=0,\tfrac12,1,\tfrac32,2,\tfrac52,\ldots$, then $f(z)=0$ identically.
\end{Corollary}

Of course, it follows directly from Carlson's theorem by substituting~$z$ with $2z$. However, much more revealing is to consider the function
\begin{gather*}
F(z)=\frac{f(z)}{\sin \pi z\,\sin\bigl(\pi\bigl[z-\tfrac12\bigr]\bigr)}
\end{gather*}
and to pass through the steps above for semicircles $\big\{\Re z\ge 0, \, |z|=n+\tfrac14\big\}$ and $\big\{\Re z\ge 0,\, |z|=n+\tfrac34\big\}$.

\begin{Corollary}\label{Corollary2} If $f(z)$ is regular and of the form ${\cal O}\big(\rme^{c|z|}\big)$, where $c<2\pi$, for $\Re z\ge 0$, and $f(z)=0$ for $z=\{n,n+a;\, n\ge 0\}$, where $0<a<1$, then $f(z)=0$ identically.
\end{Corollary}

Here we use
\begin{gather*}
F(z)=\frac{f(z)}{\sin \pi z \sin(\pi[z-a])}
\end{gather*}
and semicircles, placed in the right half-plane, with $|z|=n+\tfrac12a$ and $|z|=n+\tfrac12(a+1)$.

\begin{Corollary}\label{Corollary3} If $f(z)$ is regular and of the form ${\cal O}(\rme^{c|z|})$, where $c<3\pi$, for $\Re z\ge 0$, and $f(z)=0$ for $z=\{n,n+a,n+b;\, n\ge 0\}$, where $0<a<b<1$, then $f(z)=0$ identically.
\end{Corollary}

Here we use{\samepage
\begin{gather*}
F(z)=\frac{f(z)}{\sin \pi z \sin(\pi[z-a]) \sin(\pi[z-b])}
\end{gather*}
and semicircles with $|z|=n+\tfrac12a$, $|z|=n+\tfrac12(a+b)$, and $|z|=n+\tfrac12(b+1)$.}

The inequality $c<3\pi$ in Corollary~\ref{Corollary3} is sufficient for our purposes. Since ${\bf F}(z)-{\bf F}_0(z)$ vanishes at the points $z=\big\{n+\tfrac16, n+\tfrac12, n+\tfrac56;\, n\ge 0\big\}$, we consider the function in the half-plane $\Re z\ge\tfrac16$. The series
\begin{gather*}
{\bf F}(z)=\sum_{k=0}^\infty \frac{(z)_k\bigl(3z-\tfrac12\bigr)_k\bigl(\tfrac32-3z\bigr)_k}
{(3z)_k\bigl(\tfrac12\bigr)_k k!}\left(\frac34\right)^{k} \\
\hphantom{{\bf F}(z)}{}=\frac{\Gamma(3z)\sin\bigl(\pi\bigl[3z-\tfrac12\bigr]\bigr)}
{\sqrt{\pi} \Gamma(z)}\sum_{k=0}^\infty \frac{\Gamma(k+z)\Gamma\bigl(k+3z-\tfrac12\bigr)\Gamma\bigl(k+\tfrac32-3z\bigr)}{\Gamma(k+3z)\Gamma\bigl(k+\tfrac12\bigr)\Gamma(k+1)} \left(\frac34\right)^{k}
\end{gather*}
is absolutely convergent for any fixed $z$ in this half-plane since
\begin{gather*}
\frac{\Gamma(k+z) \Gamma\bigl(k+3z-\tfrac12\bigr)\Gamma\bigl(k+\tfrac32-3z\bigr)}
{\Gamma(k+3z)\Gamma\bigl(k+\tfrac12\bigr)\Gamma(k+1)}\sim\frac{1}{k^{2z+1/2}},\qquad k\to\infty
\end{gather*}
due to Stirling's formula \cite[equation~(5.11.7)]{DLMF}
\begin{gather*}
\Gamma(a\lambda+b)\sim\sqrt{2\pi} \rme^{-a\lambda}(a\lambda)^{a\lambda+b-1/2}, \qquad |\arg\lambda|<\pi,\qquad \lambda\to\infty, 
\end{gather*}
where $a>0$ and $b\in\Cset$ are both fixed. Thus, ${\bf F}(z)$ is holomorphic and then ${\bf F}(z)-{\bf F}_0(z)$ is holomorphic also.

The final step of our proof is to obtain an upper bound for $|{\bf F}(z)-{\bf F}_0(z)|$ as $|z|\to\infty$. Since the estimate for ${\bf F}_0(z)$ is trivial, we consider the problem for ${\bf F}(z)$ only. A major benefit of Corollary~\ref{Corollary3} use is that a very crude estimate is sufficient. Applying the integral relation \cite[{equation~(7.2.3.9)}]{Prudnikov3} which reduces ${}_{p+1}F_{q+1}$ to ${}_pF_q$, we have
\begin{gather*}
{\bf F}(z)={}_3F_2\left( \begin{matrix} z,3z-\frac12,\frac32-3z \\ 3z,\frac12\end{matrix}\,\Bigl|\,\frac34\right) \\
\hphantom{{\bf F}(z)}{}=\frac{\Gamma(3z)}{\Gamma(z)\Gamma(2z)}\int_0^1 t^{z-1}(1-t)^{2z-1}\cdot
{}_2F_1\left(\begin{matrix} 3z-\frac12,\frac32-3z \\ \frac12\end{matrix}\,\Bigl|\,\frac{3t}{4}\right) \rmd t.
\end{gather*}
Using a closed-form expression for the ${}_2F_1$-function \cite[equation~(7.3.1.90)]{Prudnikov3}
\begin{gather*}
{}_2F_1\left(\begin{matrix} a,1-a \\ \frac12\end{matrix}\,\Bigl| \,z\right)=\frac{\cos([2a-1]\arcsin\sqrt{z})}{\sqrt{1-z}},
\end{gather*}
one can obtain
\begin{gather*}
|{\bf F}(z)|\le\left|\frac{\Gamma(3z)}{\Gamma(z)\Gamma(2z)}\right| \int_0^1 t^{\Re z-1}(1-t)^{2\Re z-1}\cdot \frac{\bigl|\cos\bigl([6z-2]\arcsin\bigl(\tfrac12\sqrt{3}\,t\bigr)\bigr)\bigr|}
{\sqrt{1-\tfrac34t}}\,\rmd t.
\end{gather*}
Estimating the gamma function ratio outside the integral and the cosine term in the integrand, one finds
\begin{gather*}
\frac{\Gamma(3z)}{\Gamma(z)\Gamma(2z)}\sim\frac{\sqrt z}{\sqrt{3\pi}}\cdot \frac{3^{3z}}{2^{2z}}\quad\Rightarrow\quad
\left|\frac{\Gamma(3z)}{\Gamma(z)\Gamma(2z)}\right| ={\cal O}\left\{\sqrt{|z|}\exp\left(\Re z \ln\frac{27}{4}\right)\right\}, \quad \ln\frac{27}{4}<2,\\
|\cos(x+\rmi y)|\le\cosh y\le\rme^{|y|} \quad\Rightarrow\quad \bigl|\cos\bigl([6z-2]\arcsin\bigl(\tfrac12\sqrt{3}t\bigr)\bigr)\bigr| \le\rme^{2\pi|\Im z|}
\end{gather*}
for all $t\in[0,1]$. As result,
\begin{gather*}
|{\bf F}(z)|\le 2\rme^{2\pi|\Im z|}\cdot \left|\frac{\Gamma(3z)}{\Gamma(z)\Gamma(2z)}\right|\int_0^1 t^{\Re z-1}(1-t)^{2\Re z-1}\,\rmd t \\
\hphantom{|{\bf F}(z)|} {}=2\rme^{2\pi|\Im z|}\cdot \left|\frac{\Gamma(3z)}{\Gamma(z)\Gamma(2z)}\right| \cdot\frac{\Gamma(\Re z)\Gamma(2\Re z)}{\Gamma(3\Re z)}
\end{gather*}
and, therefore, $|{\bf F}(z)|={\cal O}\big(\rme^{c|z|}\big)$ with $c<3\pi$.

In conclusion, it is interesting to note that for proving identities (\ref{eq:HyperG_3F2_x=3/4_b})--(\ref{eq:HyperG_3F2_x=3/4_h}) the numerical part of the method above leads to the same constant, $-2$, as in (\ref{eq:Constant}), while other components of final expressions are rather
different.

\subsection*{Acknowledgements}

The authors are grateful to Professor S.K.~Suslov for his valuable comments that helped to improve the manuscript. The authors also thank the anonymous referees for their constructive criticisms and suggestions.

\pdfbookmark[1]{References}{ref}
\LastPageEnding

\end{document}